\documentclass[12pt]{article}
\usepackage{amsfonts, amssymb}

\def\AA{{\bf A}}

\def\C{\mathbb C}
\def\E{{\mathcal E}}
\def\F{{\mathcal F}}

\def\G{{\mathcal G}}

\def\H{{\mathcal H}}

\def\J{{\mathcal J}}
\def\K{{\mathcal K}}

\def\PP{{\mathbb P}}
\def\pp{{\mathfrak p}}
\def\Q{\mathbb Q}
\def\Qq{{\mathcal Q}}

\def\Rr{{\mathcal R}}
\def\U{{\mathcal U}}
\def\VV{{\bf V}}
\def\Z{\mathbb Z}

\def\Aut{\mathop{\mathfrak Aut}\nolimits}

\def\dim{\mathop{\rm dim}\nolimits}
\def\id{\mathop{\rm id}\nolimits}

\def\ker{\mathop{\rm ker}\nolimits}

\def\mod{\mathop{\rm -\!mod}\nolimits}

\def\Proj{\mathop{\bf Proj}\nolimits}
\def\proj{\mathop{\rm Proj}\nolimits}
\def\rank{\mathop{\rm rank}\nolimits}

\def\Spec{\mathop{\bf Spec}\nolimits}
\def\spec{\mathop{\rm Spec}\nolimits}
\def\Sym{\mathop{\rm Sym}\nolimits}
\def\sym{\mathop{\rm Sym}\nolimits}
\let\ov\overline
\def\im{\mathop{\rm im}\nolimits}

\def\Mor{\mathop{\mathfrak Mor}\nolimits}
\def\mor{\mathop{\rm Mor}\nolimits}
\def\qed{\hfill$\square$\bigskip}

\def\Grass{\mathop{\mathfrak Grass}\nolimits}
\def\grass{\mathop{\rm Grass}\nolimits}
\def\Hilb{\mathop{\mathfrak Hilb}\nolimits}
\def\hilb{\mathop{\rm Hilb}\nolimits}
\def\Hom{\mathop{\mathfrak Hom}\nolimits}
\def\hom{\underline{Hom}}
\def\mor{\mathop{\rm Mor}\nolimits}
\def\Quot{\mathop{\mathfrak Quot}\nolimits}
\def\quot{\mathop{\rm Quot}\nolimits}

\let\dna\downarrow
\let\hra\hookrightarrow
\let\ov\overline
\let\pa\partial
\let\un\underline

\newtheorem{theorem}{Theorem}[section]

\newtheorem{example}[theorem]{Example}
\newtheorem{lemma}[theorem]{Lemma}

\def\rem{\refstepcounter{theorem}\paragraph{Remark \thetheorem}}

\def\proof{\paragraph{Proof}}

\makeatletter \def\l@section{\@dottedtocline{1}{0em}{1.2em}} \makeatother

\begin{document}


\centerline{{\huge Construction of Hilbert and Quot Schemes}}

\bigskip

\centerline{{\Large Nitin Nitsure}}

\bigskip

{\it School of Mathematics, Tata Institute of 
Fundamental Research, Homi Bhabha Road, Mumbai 400 005, 
India. \hfill e-mail: {\tt nitsure@math.tifr.res.in}
}







\begin{abstract}
This is an expository account of 
Grothendieck's construction of Hilbert and Quot Schemes, 
following his talk 
`Techniques de construction et th\'eor\`emes d'existence
en g\'eom\'etrie alg\'ebriques IV : les sch\'emas de Hilbert',
S\'eminaire Bourbaki 221 (1960/61), 
together with further developments by Mumford and by Altman and Kleiman.
Hilbert and Quot schemes are fundamental to modern Algebraic Geometry,
in particular, for deformation theory and moduli constructions.
These notes are based on a series of six lectures in the summer 
school `Advanced Basic Algebraic Geometry', held at the 
Abdus Salam International Centre for Theoretical Physics, Trieste, 
in July 2003. 
\end{abstract}




Any scheme $X$ defines a contravariant functor $h_X$
(called the functor of points of the scheme $X$) 
from the category of schemes to the category of sets,
which associates to any scheme $T$ the set $Mor(T,X)$ of all
morphisms from $T$ to $X$. 
The scheme $X$ can be recovered (up to a unique
isomorphism) from $h_X$ by the Yoneda lemma. In fact, it is enough to know
the restriction of this functor to the full subcategory consisting of 
affine schemes, in order to recover the scheme $X$.

It is often easier to directly 
describe the functor $h_X$ than to give the scheme $X$.
Such is typically the case with various parameter schemes and
moduli schemes, or with various group-schemes over arbitrary bases,
where we can directly define a contravariant functor $F$ from 
the category of schemes to the category of sets
which would be the functor of points of the scheme in question,
without knowing in advance whether such a scheme indeed exists.

This raises the problem of representability of  
contravariant functors from 
the category of schemes to the category of sets.
An important necessary condition for 
representability come from the fact that 
the functor $h_X$ satisfies descent 
under faithfully flat quasi-compact coverings. 

(Recall that descent for a set-valued functor $F$ is the sheaf condition,
which says that if
$(f_i: U_i \to U)$ is an open cover of $U$ in the fpqc topology,
then the diagram of sets $F(U) \to \prod_i F(U_i)
\stackrel{\to}{\scriptstyle\to} \prod_{i,j}F(U_i\times_UU_j)$
is exact.)

The descent condition is often easy to verify for a given functor $F$, 
but it is not a sufficient condition for representability.

It is therefore a subtle and technically difficult problem in 
Algebraic Geometry to construct schemes which represent various
important functors, such as moduli functors.
Grothendieck addressed the issue by proving
the representability of certain basic functors, namely,
the Hilbert and Quot functors. The representing schemes
that he constructed, known as Hilbert schemes and Quot schemes, are
the foundation for proving representability of most moduli functors
(whether as schemes or as algebraic stacks).

The techniques used by Grothendieck are based on 
the theories of descent and cohomology developed by him. 
In a sequence of talks in the Bourbaki seminar, collected under the title  
`Fondements de la G\'eom\'etrie Alg\'ebriques' (see [FGA]),
he gave a sketch of the theory of descent, the 
construction of Hilbert and Quot schemes, and its application to the 
construction of Picard schemes (and also a sketch of formal schemes and
some quotient techniques).

The following notes give an expository account of 
the construction of Hilbert and Quot schemes.
We assume that the reader is familiar with the 
basics of the language of schemes and cohomology, say at the level of 
chapters 2 and 3 of Hartshorne's `Algebraic Geometry' [H]. 
Some more advanced facts about flat morphisms 
(including the local criterion for flatness) that we need are available 
in Altman and Kleiman's 
`Introduction to Grothendieck Duality Theory' [A-K 1]. 
The lecture course by Vistoli [V] on the theory of descent in 
this summer school contains in particular the background we need on descent.
Certain advanced techniques of projective geometry, 
namely Castelnuovo-Mumford regularity and flattening stratification
(to each of which we devote one lecture) are nicely given 
in Mumford's `Lectures on Curves on an Algebraic Surface' [M]. 
The book `Neron Models' by Bosch, L\"utkebohmert, Raynaud [B-L-R] 
contains a quick exposition of descent, quot schemes, and Picard schemes. 
The reader of these lecture notes is strongly urged to read Grothendieck's 
original presentation in [FGA].

\section{The Hilbert and Quot Functors}

\centerline{\large\bf The Functors $\Hilb_{\PP^n}$}

\medskip

The main problem addressed in this series of lectures, in its simplest
form, is as follows. If $S$ is a locally noetherian scheme, 
a \textbf{family of subschemes of $\PP^n$ parametrised by $S$} will mean
a closed subscheme $Y \subset \PP^n_S = \PP^n_{\Z}\times S$ 
such that $Y$ is flat over $S$. 
If $f: T\to S$ is any morphism of locally noetherian schemes, then
by pull-back we get a family $f^*(Y) = (\id \times f)^{-1}(Y) 
\subset \PP^n_T$ 
parametrised by $T$, from a family $Y$ parametrised by $S$.
This defines a contravariant functor $\Hilb_{\PP^n}$ 
from the category of all locally noetherian schemes to the category of sets, 
which associates to any $S$ the set of all such families
$$\Hilb_{\PP^n}(S) = \{ Y \subset \PP^n_S \,|\,
Y \mbox{ is flat over }S \}$$

\textbf{Question:} Is the functor $\Hilb_{\PP^n}$ representable?

Grothendieck proved that this question has an affirmative answer, that is, 
there exists a locally noetherian scheme $\hilb_{\PP^n}$ together with
a family $Z \subset \PP^n_{\Z}\times \hilb_{\PP^n}$ 
parametrised by $\hilb_{\PP^n}$,
such that any family $Y$ over $S$ is obtained as the
pull-back of $Z$ 
by a uniquely determined morphism $\varphi_Y : S \to \hilb_{\PP^n}$.
In other words, $\Hilb_{\PP^n}$ is isomorphic to 
the functor $Mor(-,\hilb_{\PP^n})$.

\medskip

\bigskip

\centerline{\large\bf The Functors $\Quot_{\oplus^r {\mathcal O}_{\PP^n}}$}
\nopagebreak
\medskip
\nopagebreak
A family $Y$ of subschemes of $\PP^n$ parametrised by $S$
is the same as 
a coherent quotient sheaf $q: {\mathcal O}_{\PP^n_S} \to {\mathcal O}_Y$ on 
$\PP^n_S$, such that ${\mathcal O}_Y$ is flat over $S$. 
This way of looking at the functor $\Hilb_{\PP^n}$
has the following fruitful generalisation.

Let $r$ be any positive integer. 
A \textbf{family of quotients of $\oplus^r {\mathcal O}_{\PP^n}$ 
parametrised by a locally noetherian scheme $S$} will mean 
a pair $(\F,q)$ consisting of 

(i) a coherent sheaf $\F$ on $\PP^n_S$ which is flat over $S$, and

(ii) a surjective ${\mathcal O}_{\PP^n_S}$-linear homomorphism of sheaves 
$q:\oplus^r {\mathcal O}_{\PP^n_S}  \to \F$.

Two such families $(\F,q)$ and $(\F,q)$ parametrised by $S$ will 
be regarded as equivalent if there exists an isomorphism $f: \F \to \F'$
which takes $q$ to $q'$, that is, the following diagram commutes.
$$\begin{array}{ccc}
\oplus^r {\mathcal O}_{\PP^n} &\stackrel{q}{\to} & \F \\
\|                  &     & ~ \dna {\scriptstyle f} \\
\oplus^r {\mathcal O}_{\PP^n} &\stackrel{q'}{\to} & \F' 
\end{array}$$
This is the same as the condition 
$\ker(q) = \ker(q')$. We will denote by $\langle \F, q \rangle$ 
an equivalence class. If $f : T\to S$ is a morphism of 
locally noetherian schemes, then pulling back the quotient
$q: \oplus^r{\mathcal O}_{\PP^n_S} \to \F$ 
under $\id\times f : \PP^n_T \to \PP^n_S$
defines a family $f^*(q) : \oplus^r {\mathcal O}_{\PP^n_T} \to f^*(\F)$ over $T$, 
which makes sense as tensor product is right-exact and preserves flatness. 
The operation of pulling back respects equivalence of families,
therefore it gives rise to a contravariant functor
$\Quot_{\oplus^r {\mathcal O}_{\PP^n}}$ from the category of
all locally noetherian schemes to the category of sets, 
by putting 
$$\Quot_{\oplus^r {\mathcal O}_{\PP^n}}(S) = 
\{\mbox{ All } \langle \F, q \rangle \mbox{ parametrised by }S\}$$

It is immediate that the functor $\Quot_{\oplus^r {\mathcal O}_{\PP^n}}$
satisfies faithfully flat descent. It was proved by
Grothendieck that in fact the above functor
is representable on the category of all locally noetherian schemes 
by a scheme $\quot_{\oplus^r{\mathcal O}_{\PP^n}}$.

\medskip

\bigskip

\centerline{\large\bf The Functors $\Hilb_{X/S}$ and $\Quot_{E/X/S}$}

\medskip

The above functors $\Hilb_{\PP^n}$ and $\Quot_{\oplus^r {\mathcal O}_{\PP^n}}$
admit the following simple generalisations. 
Let $S$ be a noetherian scheme and let $X\to S$ be a finite type
scheme over it. Let $E$ be a coherent sheaf on $X$. 
Let $Sch_S$ denote the category of all locally noetherian schemes
over $S$. For any $T\to S$ in $Sch_S$, a \textbf{family of quotients of $E$ 
parametrised by $T$} will mean 
a pair $(\F,q)$ consisting of 

(i) a coherent sheaf $\F$ on $X_T = X\times_ST$ 
such that the schematic support of $\F$ is proper over $T$ 
and $\F$ is flat over $T$, together with 

(ii) a surjective ${\mathcal O}_{X_T}$-linear homomorphism of sheaves 
$q:E_T  \to \F$ where $E_T$ is the pull-back of $E$ under
the projection $X_T \to X$.

Two such families $(\F,q)$ and $(\F,q)$ parametrised by $T$ will 
be regarded as equivalent if $\ker(q) = \ker(q')$), 
and $\langle \F, q \rangle$ will denote 
an equivalence class. Then as properness and flatness are preserved
by base-change, and as tensor-product is right exact, 
the pull-back of $\langle \F, q \rangle$ under an $S$-morphism $T'\to T$ is
well-defined, which gives
a set-valued contravariant functor 
$\Quot_{E/X/S} : Sch_S \to Sets$ under which   
$$T\mapsto  \{\mbox{ All } \langle \F, q \rangle \mbox{ parametrised by }T\}$$

When $E = {\mathcal O}_X$, the functor $\Quot_{{\mathcal O}_X/X/S} : Sch_S \to Sets$ 
associates to $T$ the set of all closed subschemes $Y\subset X_T$
that are proper and flat over $T$. We denote this functor by $\Hilb_{X/S}$.

Note in particular that we have
$$\Hilb_{\PP^n} = \Hilb_{\PP^n_{\Z}/\spec \Z} \mbox{ and }
\Quot_{\oplus^r {\mathcal O}_{\PP^n}} = 
\Quot_{\oplus^r {\mathcal O}_{\PP^n_{\Z}}/\PP^n_{\Z}/\spec \Z}$$

It is clear that the functors $\Quot_{E/X/S}$ and $\Hilb_{X/S}$
satisfy faithfully flat descent, 
so it makes sense to pose the question of their
representability.

\medskip

\bigskip

\centerline{\large\bf Stratification by Hilbert Polynomials}

\medskip

Let $X$ be a finite type scheme over
a field $k$, together with a line bundle $L$.
Recall that if $F$ is a coherent sheaf on 
$X$ whose support is proper over $k$, then
the \textbf{Hilbert polynomial} ${\Phi}\in \Q[\lambda]$ of $F$ is defined by
the function 
$${\Phi}(m) = \chi(F(m)) = 
\sum_{i=0}^n (-1)^i \dim_k H^i(X, F\otimes L^{\otimes m})$$
where the dimensions of the cohomologies are finite because of the 
coherence and properness conditions. The fact that $\chi(F(m))$
is indeed a polynomial in $m$ under the above assumption is a special
case of what is known as Snapper's Lemma (see Kleiman [K] for a proof).

Let $X\to S$ be a finite type morphism of noetherian schemes,
and let $L$ be a line bundle on $X$. Let $\F$
be any coherent sheaf on 
$X$ whose schematic support is proper over $S$. Then for each $s\in S$,
we get a polynomial ${\Phi}_s \in  \Q[\lambda]$ 
which is the Hilbert polynomial of
the restriction $\F_s = F|_{X_s}$ of $\F$ to the fiber $X_s$ over $s$,
calculated with respect to the line bundle $L_s = L|_{X_s}$. 
If $\F$ is flat over $S$ then the function $s\mapsto {\Phi}_s$ 
from the set of points of $S$ to the polynomial ring $\Q[\lambda]$ is 
known to be locally constant on $S$.

This shows that the functor $\Quot_{E/X/S}$ naturally decomposes
as a co-product 
$$\Quot_{E/X/S} = \coprod_{{\Phi}\in \Q[\lambda]}\, \Quot_{E/X/S}^{\Phi,L}$$
where for any polynomial ${\Phi}\in \Q[\lambda]$, the functor 
$\Quot_{E/X/S}^{\Phi,L}$ associates
to any $T$ the set of all equivalence classes of 
families $\langle \F,q \rangle$ such that at each $t\in T$ the 
Hilbert polynomial of the restriction $\F_t$, calculated
using the pull-back of $L$, is ${\Phi}$. Correspondingly,
the representing scheme  $\quot_{E/X/S}$, when it exists, naturally decomposes
as a co-product
$$\quot_{E/X/S} = \coprod_{{\Phi}\in \Q[\lambda]}\, \quot_{E/X/S}^{\Phi,L}$$

\textbf{Note } 
We will generally take $X$ to be (quasi-)projective over $S$,
and $L$ to be a relatively very ample line bundle. Then indeed 
the Hilbert and Quot functors are representable by schemes,
but not in general.

\medskip

\bigskip

\centerline{\large\bf Elementary Examples, Exercises}

\medskip

\textbf{(1) } {\bf $\PP^n_{\Z}$ as a Quot scheme } 
Show that the scheme $\PP^n_{\Z} = \proj \Z[x_0,\ldots,x_n]$ represents the 
functor $\varphi$ from schemes to sets, which associates
to any $S$ the set of all equivalence classes $\langle \F, q\rangle$ 
of quotients
$q: \oplus^{n+1}{\mathcal O}_S \to \F$, where $\F$ is an invertible ${\mathcal O}_S$-module. 
As coherent sheaves on $S$ which are ${\mathcal O}_S$-flat with each
fiber $1$-dimensional are exactly the locally free sheaves on $S$ 
of rank $1$, it follows that $\varphi$ is the functor 
$\Quot_{\oplus^{n+1}{\mathcal O}_{\Z}/\Z/\Z}^{1, {\mathcal O}_{\Z}}$
(where in some places we write just $\Z$ for $\spec \Z$ for simplicity). 
This shows that
$\quot_{\oplus^{n+1}{\mathcal O}_{\Z}/\Z/\Z}^{1, {\mathcal O}_{\Z}} 
= \PP^n_{\Z}$. 
Under this identification, show that the universal family on
$\quot_{\oplus^{n+1}{\mathcal O}_{\Z}/\Z/\Z}^{1, {\mathcal O}_{\Z}}$ 
is the tautological quotient 
$\oplus^{n+1}{\mathcal O}_{\PP^n_{\Z}} \to {\mathcal O}_{\PP^n_{\Z}}(1)$
 
More generally, show that if $E$ is a locally free sheaf on a 
noetherian scheme $S$, the functor $\Quot_{E/S/S}^{1, {\mathcal O}_S}$ is represented 
by the $S$-scheme $\PP(E) = \Proj \sym_{{\mathcal O}_S} E$, 
with the tautological quotient 
$\pi^*(E) \to {\mathcal O}_{\PP(E)}(1)$ as the universal family.

\medskip

\textbf{(2) } \textbf{Grassmannian as a Quot scheme } 
For any integers $r\ge d\ge 1$, an explicit construction
the Grassmannian scheme
$\grass(r,d)$ over $\Z$, together with the tautological 
quotient $u:\oplus^r\,{\mathcal O}_{\grass(r,d)} \to \U$ where \
$\U$ is a rank $d$ locally free sheaf on $\grass(r,d)$, 
has been given at the end of this section. A proof of the 
properness of $\pi : \grass(r,d) \to \spec \Z$ is given there, 
together with a closed embedding $\grass(r,d)\hra 
\PP(\pi_*\det\U) = \PP^m_{\Z}$ where $m={r \choose d} - 1$.

Show that $\grass(r,d)$ together with the quotient 
$u:\oplus^r\,{\mathcal O}_{\grass(r,d)} \to \U$ represents the 
contravariant functor 
$$\Grass(r,d) = \Quot_{\oplus^r\,{\mathcal O}_{\Z}/\Z/\Z}^{d, {\mathcal O}_{\Z}}$$ 
from schemes to sets, which associates to any $T$ the
set of all equivalence classes 
$\langle \F, q\rangle$ of quotients $q: \oplus^r\,{\mathcal O}_T \to \F$
where $\F$ is a locally free sheaf on $T$ of rank $d$.
Therefore, $\quot_{\oplus^r\,{\mathcal O}_{\Z}/\Z/\Z}^{d, {\mathcal O}_{\Z}}$ 
exists, and equals $\grass(r,d)$.

\textbf{Grassmannian of a vector bundle} 
Show that for any ring
$A$, the action of the group $GL_r(A)$ on the free module $\oplus^r A$ 
induces an action of $GL_r(A)$ on 
the set $\Grass(r,d)(A)$, such that for any ring homomorphism $A\to B$,
the set-map $\Grass(r,d)(A)\to\Grass(r,d)(B)$ is equivariant with respect
to the group homomorphism $GL_r(A) \to GL_r(B)$. (In schematic terms,
this means we have an action of the group-scheme $GL_{r,\Z}$ on 
$\grass(r,d)$.) 

Using the above show that, more generally, if $S$ is a scheme and 
$E$ is a locally free ${\mathcal O}_S$-module of rank $r$, 
the functor $\Grass(E,d)=\Quot_{E/S/S}^{d,{\mathcal O}_S}$ 
on all $S$-schemes which 
by definition associates to any $T$ 
the set of all equivalence classes $\langle \F, q\rangle$ of quotients  
$q: E_T \to \F$ where $\F$ is a locally free sheaf on $T$ of rank $d$,
is representable. The representing scheme
is denoted by $\grass(E,d)$ and is called the rank $d$ 
relative Grassmannian of $E$ over $S$.
It parametrises a universal quotient $\pi^*E \to \F$
where $\pi : \grass(E,d) \to S$ is the projection.
Show that the determinant line bundle $\bigwedge^d \F$ on
$\grass(E,d)$ is relatively very ample over $S$, 
and it gives a closed embedding
$\grass(E,d) \hra \PP(\pi_* \bigwedge^d \F) \subset \PP(\bigwedge^d E)$. 
(The properness of the embedding follows from the 
properness of $\pi : \grass(E,d) \to S$, which follows locally
over $S$ by base-change from  properness of 
$\grass(r,d)$ over $\Z$ -- see Exercise \textbf{(5)}
or \textbf{(7)} below.)

\textbf{Grassmannian of a coherent sheaf} 
If $E$ is a coherent sheaf on $S$, not necessarily locally free, 
then by definition the functor $\Grass(E,d)=\Quot_{E/S/S}^{d,{\mathcal O}_S}$ 
on all $S$-schemes associates to any $T$ 
the set of all equivalence classes $\langle \F, q\rangle$ of quotients  
$q: E_T \to \F$ where $\F$ is a locally free on $T$ of rank $d$.
If $r: E' \to E$ is a surjection of coherent sheaves on $S$, then show that
the induced morphism of functors $\Grass(E,d)\to \Grass(E',d)$,
which sends $\langle \F, q\rangle \mapsto \langle \F, q\circ r\rangle$,
is a closed embedding. 
From this, by locally expressing a coherent
sheaf as a quotient of a vector bundle, show that $\Grass(E,d)$ 
is representable even when $E$ is a coherent sheaf on $S$ which is not 
necessarily locally free. The representing scheme $\grass(E,d)$ is proper
over $S$, as locally over $S$ it is a closed
subscheme of the Grassmannian of a vector bundle. 
Show by arguing locally over $S$ that the line bundle
$\bigwedge^d \F$ on $\grass(E,d)$ is relatively very ample over $S$, 
and therefore by using properness 
conclude that $\grass(E,d)$ is projective over $S$.

\medskip

\textbf{(3) } \textbf{Grassmannian as a Hilbert scheme }
Let $\Phi = 1 \in \Q[\lambda]$. Then the Hilbert scheme
$\hilb_{\PP^n}^{1,{\mathcal O}(1)}$ is $\PP^n_{\Z}$ itself.
More generally, let $\Phi_r = {r + \lambda \choose r}\in \Q[\lambda]$
where $r \ge 0$. 
The Hilbert scheme $\hilb_{\PP^n}^{\Phi_r,{\mathcal O}(1)}$ is isomorphic to the 
Grassmannian scheme $\grass(n+1, r+1)$ over $\Z$. 
This can be seen via the following steps, whose detailed verification 
is left to the reader as an exercise.

(i) The Grassmannian scheme 
$\grass(n+1, r+1)$ over $\Z$ parametrises a tautological family of 
subschemes of $\PP^n$ with Hilbert polynomial $\Phi_r$. 
Therefore we get a natural transformation 
$h_{\grass(n+1, r+1)} \to \Hilb_{\PP^n}^{\Phi_r,{\mathcal O}(1)}$.

(ii) Any closed subscheme $Y \subset \PP^n_k$
with Hilbert polynomial ${\Phi}_r$, where $k$ is any field, 
is isomorphic to $\PP^r_k$ embedded linearly in $\PP^n_k$ over $k$.
If $V$ is a vector bundle over a noetherian base $S$, and if
$Y\subset \PP(V)$ is a closed subscheme flat over $S$ with each
schematic fiber $Y_s$ an $r$-dimensional 
linear subspace of the projective space
$\PP(V_s)$, then $Y$ is defines
a rank $r+1$ quotient vector bundle 
$V = \pi_*{\mathcal O}_{\PP(V)}(1) \to \pi_*{\mathcal O}_Y(1)$  
where $\pi : \PP(V) \to S$ denotes the projection. 
This gives a natural transformation
$\Hilb_{\PP^n}^{{\Phi}_r,{\mathcal O}(1)} \to h_{\grass(n+1, r+1)}$.

(iii) The above two natural transformations are inverses of each other.

\medskip

\textbf{(4) } \textbf{Hilbert scheme of hypersurfaces in $\PP^n$ }
Let  
${\Phi}_d = {n+\lambda \choose n} - {n -d + \lambda\choose n}
\in \Q[\lambda]$
where $d\ge 1$. 
The Hilbert scheme $\hilb_{\PP^n}^{{\Phi}_d,{\mathcal O}(1)}$ is isomorphic to 
$\PP^m_{\Z}$ where $m = {n+d \choose d} - 1$. 
This can be seen from the following steps, which are
left as exercises.

(i) Any closed subscheme $Y \subset \PP^n_k$
with Hilbert polynomial ${\Phi}_d$, where $k$ is any field, 
is a hypersurface of degree $d$ in $\PP^n_k$. \textit{Hint:}
If $Y \subset \PP^n_k$ is a closed subscheme with
Hilbert polynomial of degree $n-1$, then show that
the schematic closure $Z$ of the hight $1$
primary components is a hypersurface in $\PP^n_k$ with
$\deg(Z) = \deg(Y)$. 

(ii) Any family $Y \subset \PP^n_S$ is a Cartier divisor
in $\PP^n_S$. It gives rise to a line subbundle 
$\pi_* (I_Y \otimes {\mathcal O}_{\PP^n_S}(d)) \subset \pi_* {\mathcal O}_{\PP^n_S}(d)$,
which defines a natural morphism 
$f_Y : S \to \PP^m_{\Z}$ where $m = {n+d \choose  d} - 1$.
This gives a morphism of functors 
$\Hilb_{\PP^n}^{{\Phi}_d}\to \PP^m$ where we denote
$h_{\PP^m_{\Z}}$ simply by $\PP^m$.

(iii) The scheme $\PP^m_{\Z}$ parametrises a tautological family of 
hypersurfaces of degree $d$, which gives a 
morphism of functors 
$\PP^m \to \Hilb_{\PP^n}^{{\Phi}_d,{\mathcal O}(1)}$ in the reverse direction.
These are inverses of each other.

\medskip

\textbf{(5) } \textbf{Base-change property of Hilbert and Quot schemes} 
Let $S$ be a noetherian scheme, $X$ a finite-type scheme over
$S$, and $E$ a coherent sheaf on $X$. 
If $T\to S$ is
a morphism of noetherian schemes, then show that there is a natural
isomorphism of functors
$\Quot_{E_T/X_T/T} \to \Quot_{E/X/S}\times_{h_S} h_T$.
Consequently, if $\quot_{E/X/S}$ exists, then so does
$\quot_{E_T/X_T/T}$, which is naturally isomorphic to 
$\quot_{E/X/S}\times_ST$. One can prove a similar
statement involving $\Quot_{E/X/S}^{\Phi,L}$.
In particular, $\hilb_{X/S}$ and $\hilb_{X/S}^{\Phi,L}$, 
when they exist, base-change correctly.

\medskip

\textbf{(6) } \textbf{Descent condition in the fpqc topology}
If $U$ is an $S$-scheme and $(f_i: U_i \to U)$ is an open cover of $U$ 
in the fpqc topology, then show that the following 
sequence of sets is exact:
$$\Quot_{E/X/S}(U) \to \prod_i \Quot_{E/X/S}(U_i)
\stackrel{\to}{\scriptstyle\to} \prod_{i,j}\Quot_{E/X/S}(U_i\times_UU_j)$$

\medskip

\textbf{(7) } \textbf{Valuative criterion for properness }
When $X\to S$ is proper, show that the morphism of functors
$\Quot_{E/X/S} \to h_S$ satisfies the valuative criterion of properness
with respect to discrete valuation rings, that is, 
if $R$ is a discrete valuation ring together
with a given morphism $\spec R \to S$ making it an $S$-scheme, 
show that the restriction map 
$\Quot_{E/X/S}(\spec R) \to \Quot_{E/X/S}(\spec K)$
is bijective, where $K$ is the quotient field of $R$ 
and $\spec K$ is regarded as an $S$-scheme in the obvious way.

\medskip

\textbf{(8) } \textbf{Counterexample of Hironaka }
Hironaka constructed a $3$-dimensional smooth proper scheme $X$ over 
complex numbers $\C$, together with a free action of the group
$G = \Z/(2)$, for which the quotient $X/G$ does not exist as a scheme.
(See Example 3.4.1 in Hartshorne [H] Appendix B for construction of
$X$. We leave the definition of the $G$ action and the proof that
$X/G$ does not exist to the reader.)
In particular, this means the Hilbert functor $\Hilb_{X/\C}$
is not representable by a scheme.

\medskip

\bigskip


\centerline{\large\bf Construction of Grassmannian}

\medskip

The following explicit construction of 
the Grassmannian scheme $\grass(r,d)$ over $\Z$ is best understood 
as the construction of a quotient
$GL_{d,\Z}\backslash V$, where $V$ is the scheme of all
$d\times r$ matrices of rank $d$, and the group-scheme
$GL_{d,\Z}$ acts on $V$ on the left by matrix multiplication. However,
we will not use the language of group-scheme actions 
here, instead, we give a direct elementary construction of 
the Grassmannian scheme.

The reader can take $d=1$ in what follows, in a first reading, to
get the special case $\grass(r,1) = \PP^{r-1}_{\Z}$, which 
has another construction as $\proj \Z[x_1,\ldots,x_r]$.

\medskip

\textbf{Construction by gluing together affine patches }
For any integers $r\ge d\ge 1$, the Grassmannian scheme
$\grass(r,d)$ over $\Z$, together with the tautological 
quotient $u:\oplus^r\,{\mathcal O}_{\grass(r,d)} \to \U$ where \
$\U$ is a rank $d$ locally free sheaf on $\grass(r,d)$, 
can be explicitly constructed as follows.

If $M$ is a $d\times r$-matrix, and 
$I\subset \{1,\ldots, r\}$ with cardinality $\#(I)$ equal 
to $d$, the $I$ th minor $M_I$ of $M$ will 
mean the $d\times d$ minor of $M$ whose columns are indexed by
$I$.

For any subset $I\subset \{1,\ldots, r\}$ with $\#(I) =d$, 
consider the $d\times r$ matrix $X^I$ whose $I$ the minor
$X^I_I$ is the $d\times d$ identity matrix $1_{d\times d}$, while the remaining
entries of $X^I$ are independent variables $x^I_{p,q}$ over $\Z$. 
Let $\Z[X^I]$ denote the polynomial ring in the variables
$x^I_{p,q}$, and let $U^I = \spec \Z[X^I]$, which is non-canonically 
isomorphic to the affine space $\AA^{d(r-d)}_{\Z}$. 

For any $J\subset \{1,\ldots, r\}$ with $\#(J) =d$, let 
$P^I_J = \det(X^I_J) \in \Z[X^I]$
where $X^I_J$ is the $J$ th minor of $X^I$.
Let $U^I_J = \spec \Z[X^I, 1/P^I_J]$
the open subscheme of $U^I$ where $P^I_J$ is invertible. This means
the $d\times d$-matrix $X^I_J$ admits an inverse $(X^I_J)^{-1}$
on $U^I_J$. 

For any $I$ and $J$, a ring homomorphism 
$\theta_{I,J} : \Z[X^J, 1/P^J_I] \to \Z[X^I, 1/P^I_J]$ is
defined as follows. The images of the variables 
$x^J_{p,q}$ are given by the entries of the matrix formula
$\theta_{I,J}(X^J) = (X^I_J) ^{-1}\,X^I$.
In particular, we have $\theta_{I,J}(P^J_I) = 1/P^I_J$,
so the map extends to $\Z[X^J, 1/P^J_I]$.

Note that $\theta_{I,I}$ is identity on $U^I_I = U^I$, and 
we leave it to the reader to verify that 
for any three subsets $I$, $J$ and $K$ of $\{ 1,\ldots, r\}$ of 
cardinality $d$, the co-cycle condition 
$\theta_{I,K}= \theta_{I,J}\theta_{J,K}$
is satisfied. Therefore the schemes 
$U^I$, as $I$ varies over all the 
${r\choose d}$ different 
subsets of  $\{ 1,\ldots, r\}$ of cardinality $d$,
can be glued together by the co-cycle $(\theta_{I,J})$ to form a
finite-type scheme $\grass(r,d)$ over $\Z$. As each $U^I$ is
isomorphic to $\AA^{d(r-d)}_{\Z}$, it follows that
$\grass(r,d) \to \spec \Z$ is smooth of relative dimension $d(r-d)$.

\medskip

\textbf{Separatedness }
The intersection of the diagonal of $\grass(r,d)$ with $U^I\times U^J$
can be seen to be the closed subscheme $\Delta_{I,J}\subset U^I\times U^J$ 
defined by entries of the matrix formula $X^J_I\,X^I - X^J =0$, 
and so $\grass(r,d)$ is a separated scheme.

\medskip

\textbf{Properness }
We now show that $\pi : \grass(r,d) \to \spec \Z$ is proper. It is enough to
verify the valuative criterion of properness for discrete valuation rings.
Let $\Rr$ be a dvr, $\K$ its quotient field, and let
$\varphi : \spec \K \to \grass(r,d)$ be a morphism. This is given by 
a ring homomorphism $f : \Z[X^I] \to \K$ for some $I$. Having fixed
one such $I$, next choose $J$ such that $\nu(f(P^I_J))$ is minimum,
where $\nu : \K \to \Z\bigcup \{ \infty\}$ denotes 
the discrete valuation. As $P^I_I = 1$, note that
$\nu(f(P^I_J)) \le 0$, therefore
$f(P^I_J) \ne 0$ in $\K$ and so the matrix $f(X^I_J)$
lies in $GL_d(\K)$.

Now consider the homomorphism
$g : \Z[X^J] \to \K$ defined by entries of the matrix formula
$$g(X^J) = f((X^I_J)^{-1}\,X^I)$$
Then $g$ defines the same morphism $\varphi : \spec \K \to \grass(r,d)$,
and moreover all $d\times d$ minors $X^J_K$ satisfy
$\nu(g(P^J_K)) \ge 0$. 
As the minor $X^J_J$ is identity, it follows from the above that in fact
$\nu(g(x^J_{p,q}))\ge 0$ for all entries of $X^J$. 
Therefore, the map $g : \Z[X^J] \to \K$
factors uniquely via $\Rr \subset \K$. The resulting morphism of
schemes $\spec \Rr \to U^J \hra \grass(r,d)$ prolongs
$\varphi : \spec \K \to \grass(r,d)$. We have already checked
separatedness of $\grass(r,d)$, so now we see that 
$\grass(r,d) \to \spec \Z$ is proper.

\medskip

\textbf{Universal quotient }
We next define a rank $d$ locally free sheaf $\U$ on $\grass(r,d)$
together with a surjective homomorphism $\oplus^r\,{\mathcal O}_{\grass(r,d)} \to \U$.
On each $U^I$ we define a surjective homomorphism 
$u^I :\oplus^r\,{\mathcal O}_{U^I} \to  \oplus^d\,{\mathcal O}_{U^I}$
by the matrix $X^I$. Compatible with the co-cycle $(\theta_{I,J})$
for gluing the affine pieces $U^I$, we give gluing data
$(g_{I,J})$ for gluing together the trivial bundles $\oplus^d\,{\mathcal O}_{U^I}$
by putting 
$$g_{I,J} = (X^I_J)^{-1} \in GL_d(U^I_J)$$
This is compatible with the homomorphisms $u^I$, 
so we get a surjective homomorphism 
$u : \oplus^r\,{\mathcal O}_{\grass(r,d)} \to \U$.

\medskip
\textbf{Projective embedding }
As $\U$ is given by the transition functions $g_{I,J}$ described above, 
the determinant line bundle 
$\det(\U)$ is given by the transition functions
$\det(g_{I,J}) = 1/P^I_J \in GL_1(U^I_J)$.
For each $I$, we define a global section 
$$\sigma_I \in \Gamma(\grass(r,d),\det(\U))$$ 
by putting $\sigma_I|_{U^J} = P^J_I \in 
\Gamma(U^J,{\mathcal O}_{U^J})$ in terms of the trivialization over the open cover
$(U^J)$. 
We leave it to the reader to verify that the sections $\sigma_I$ form
a linear system which is base point free and separates points 
relative to $\spec \Z$, and so gives an embedding of 
$\grass(r,d)$ into $\PP^m_{\Z}$ where $m = {r\choose d}-1$.
This is a closed embedding by the properness of 
$\pi:\grass(r,d) \to \spec \Z$. In particular, $\det(\U)$
is a relatively very ample line bundle on $\grass(r,d)$ over $\Z$.

\medskip

{\footnotesize
\textbf{Note } The $\sigma_I$ are known as the Pl\"ucker coordinates, 
and these satisfy certain quadratic polynomials known as the Pl\"ucker 
relations, which define the projective image of the Grassmannian.
We will not need these facts. }

\section{Castelnuovo-Mumford Regularity}

Mumford's deployment of $m$-regularity led to a simplification in the 
construction of Quot schemes. The original construction of Grothendieck 
had instead relied on Chow coordinates.

Let $k$ be a field and let $\F$ be a coherent sheaf on the 
projective space $\PP^n$ over $k$. Let $m$ be an integer.
The sheaf $\F$ is said to be \textbf{$m$-regular} if we have
$$H^i(\PP^n, \F(m-i)) =0\mbox{ for each }i \ge 1.$$
The definition, which may look strange at first sight, is suitable
for making inductive arguments on $n=\dim(\PP^n)$
by restriction to a suitable hyperplane. If $H\subset \PP^n$ is
a hyperplane which does not contain any associated point of $\F$,
then we have a short exact sheaf sequence
$$0 \to \F(m-i-1) \stackrel{\alpha}{\to} \F(m-i) \to \F_H(m-i)\to 0$$
where the map $\alpha$ is locally given by multiplication with
a defining equation of $H$, hence is injective.
The resulting long exact cohomology sequence
$$\ldots \to H^i(\PP^n, \F(m-i)) \to H^i(\PP^n, \F_H(m-i)) 
\to H^{i+1}(\PP^n, \F(m-i-1)) \to \ldots$$
shows that if $\F$ is $m$-regular, then so is its restriction $\F_H$
(with the \textit{same value\,} for $m$) 
to a hyperplane $H\simeq \PP^{n-1}$ which does not contain any 
associated point of $\F$.
Note that whenever $\F$ is coherent, the set of associated points of
$\F$ is finite, so there will exist at least one such hyperplane $H$
when the field $k$ is infinite.

The following lemma is due to Castelnuovo, according to 
Mumford's {\sl Curves on a surface}.

\begin{lemma}\label{Castelnuovo} 
If $\F$ is an $m$-regular sheaf on $\PP^n$ then
the following statements hold:

\textbf{(a)} The canonical map $H^0(\PP^n,{\mathcal O}_{\PP^n}(1))\otimes 
H^0(\PP^n,\F(r)) \to H^0(\PP^n,\F(r+1))$ is surjective
whenever $r\ge m$.

\textbf{(b)} We have $H^i(\PP^n, \F(r))=0$ whenever $i \ge 1$ and $r\ge m-i$.
In other words, if $\F$ is $m$-regular, then it is $m'$-regular for all
$m' \ge m$.

\textbf{(c)} The sheaf $\F(r)$ is generated by its global sections, 
and all its higher cohomologies vanish, whenever $r\ge m$.
\end{lemma}

\proof As the cohomologies base-change correctly under a field extension, 
we can assume that the field $k$ is infinite.
We argue by induction on $n$. The statements
\textbf{(a)}, \textbf{(b)} and \textbf{(c)} clearly hold when $n=0$, so
next let $n\ge 1$. 
As $k$ is infinite, there exists a hyperplane $H$ 
which does not contain any associated point of $\F$, so that
the restriction $\F_H$ is again $m$-regular as explained above.
As $H$ is isomorphic to $\PP^{n-1}_k$, by the 
inductive hypothesis the assertions of the lemma hold for the sheaf $\F_H$.

When $r=m-i$, the equality $H^i(\PP^n, \F(r))=0$ in statement \textbf{(b)}
follows for all $n\ge 0$ by definition of $m$-regularity. 
To prove \textbf{(b)}, we now proceed by induction on $r$ where 
$r\ge m-i+1$. Consider the exact sequence 
$$H^i(\PP^n, \F(r-1)) \to H^i(\PP^n, \F(r)) \to 
H^i(H, \F_H(r))$$
By inductive hypothesis for $r-1$ 
the first term is zero, 
while by inductive hypothesis for $n-1$
the last term is zero, which shows that the middle term is zero, 
completing the proof of \textbf{(b)}. 

Now consider the commutative diagram
{\footnotesize
$$\begin{array}{lccc}
& H^0(\PP^n, \F(r))\otimes H^0(\PP^n,{\mathcal O}_{\PP^n}(1))& 
\stackrel{\sigma}{\to}& H^0(H, \F_H(r))\otimes H^0(H, {\mathcal O}_H(1)) \\
& ~\dna {\scriptstyle \mu}& & ~\dna {\scriptstyle \tau}\\
H^0(\PP^n, \F(r))~~~~ \stackrel{\alpha}{\to} & H^0(\PP^n, \F(r+1))
& \stackrel{\nu_{r+1}}{\to} & H^0(H, \F_H(r+1))
\end{array}$$
}
The top map $\sigma$ is surjective, for the following reason: 
By $m$-regularity of $\F$ and using the statement \textbf{(b)}
already proved, we see that 
$H^1(\PP^n, \F(r-1))=0$ for $r\ge m$, and so the restriction map
$\nu_r: H^0(\PP^n, \F(r)) \to H^0(H, \F_H(r))$ is surjective. 
Also, the restriction map
$\rho : H^0(\PP^n,{\mathcal O}_{\PP^n}(1)) \to H^0(H, {\mathcal O}_H(1))$ is surjective. 
Therefore the tensor product $\sigma = \nu_r \otimes \rho$ 
of these two maps is surjective. 

The second vertical map $\tau$ is surjective by inductive
hypothesis for $n-1 = \dim(H)$.

Therefore, the composite $\tau\circ \sigma$ is surjective,
so the composite $\nu_{r+1}\circ\mu$ is surjective,
hence $H^0(\PP^n, \F(r+1)) = \im(\mu) + \ker(\nu_{r+1})$.
As the bottom row is exact, we get 
$H^0(\PP^n, \F(r+1)) =  \im(\mu) + \im(\alpha)$. 
However, we have $\im(\alpha) \subset \im(\mu)$, as the map
$\alpha$ is given by tensoring with a certain section of 
${\mathcal O}_{\PP^n}(1)$ (which has divisor $H$). 
Therefore, $H^0(\PP^n, \F(r+1)) = \im(\mu)$.
This completes the proof of \textbf{(a)} for all $n$.

To prove \textbf{(c)}, consider the map $H^0(\PP^n, \F(r))
\otimes H^0(\PP^n,{\mathcal O}_{\PP^n}(p)) \to
H^0(\PP^n, \F(r+p))$, 
which is
surjective for $r\ge m$ and $p\ge 0$ as follows from 
a repeated use of \textbf{(a)}. 
For $p\gg 0$, we know that 
$H^0(\PP^n, \F(r+p))$ is generated by its global sections.
It follows that $H^0(\PP^n, \F(r))$ is also generated 
by its global sections for $r\ge m$. We already know from \textbf{(b)}
that $H^i(\PP^n, \F(r))=0$ for $i\ge 1$ when $r\ge m$.
This proves \textbf{(c)}, completing the proof of the lemma.
\qed

\rem\label{surjectivity of restriction} 
The following fact, based on the
diagram used in the course of the above proof, 
will be useful later: With notation as above, 
let the restriction map
$\nu_r: H^0(\PP^n, \F(r)) \to 
H^0(H, \F_H(r))$ be surjective. Also, let $\F_H$ be $r$-regular, so that
by Lemma \ref{Castelnuovo}.\textbf{(a)} the map 
$H^0(H, {\mathcal O}_H(1))\otimes H^0(H, \F_H(r))\to H^0(H, \F_H(r+1))$ 
is surjective. Then the restriction map 
$\nu_{r+1} : H^0(\PP^n, \F(r+1)) \to  H^0(H, \F_H(r+1))$
is again surjective. As a consequence, if
$\F_H$ is $m$ regular and if for some $r \ge m$ the 
restriction map $\nu_r : H^0(\PP^n, \F(r)) \to 
H^0(H, \F_H(r))$ is surjective, then the restriction map 
$\nu_p : H^0(\PP^n, \F(p)) \to 
H^0(H, \F_H(p))$ is surjective for all $p\ge r$.

\medskip

\textbf{Exercise } Find all the values of $m$ for which the 
invertible sheaf ${\mathcal O}_{\PP^n}(r)$ is $m$-regular.

\medskip

\textbf{Exercise } Suppose $0\to \F' \to \F \to \F'' \to 0$
is an exact sequence of coherent sheaves on $\PP^n$.
Show that if $\F'$ and $\F''$ are $m$-regular, then $\F$ is also
$m$-regular, if $\F'$ is $(m+1)$-regular and $\F$ is $m$-regular, 
then $\F''$ is $m$-regular, and 
if $\F$ is $m$-regular and $\F''$ is $(m-1)$-regular, 
then $\F'$ is $m$-regular.

\bigskip

The use of $m$-regularity for making Quot schemes is via the 
following theorem. 

\begin{theorem}\label{Mumford} {\rm (Mumford)} 
For any non-negative integers $p$ and $n$,
there exists a polynomial $F_{p,n}$ in $n+1$ variables
with integral coefficients, which has the following property:

Let $k$ be any field, and let $\PP^n$ denote the $n$-dimensional 
projective space over $k$. Let
$\F$ be any coherent sheaf 
on $\PP^n$, which is isomorphic to a
subsheaf of $\oplus ^p{\mathcal O}_{\PP^n}$.
Let the Hilbert polynomial of $\F$ be written in terms of binomial 
coefficients as
$$\chi(\F(r)) = \sum_{i=0}^n a_i \, {r\choose i}$$
where $a_0,\ldots,a_n \in \Z$.

Then $\F$ is $m$-regular, where 
$m=F_{p,n}(a_0,\ldots,a_n)$.
\end{theorem}

\proof (Following Mumford [M])
As before, we can assume that $k$ is infinite. We argue by induction on $n$. 
When $n=0$, clearly we can take $F_{p,0}$ to be any polynomial.
Next, let $n\ge 1$. 
Let $H\subset \PP^n$ be a hyperplane which does not contain
any of the finitely many associated points of
$\oplus^p{\mathcal O}_{\PP^n}/\F$ (such an $H$ exists as $k$ is infinite).
Then the following torsion sheaf vanishes:
$$\un{Tor}^1_{{\mathcal O}_{\PP^n}}({\mathcal O}_H,\, \oplus ^p{\mathcal O}_{\PP^n}/\F)=0$$ 
Therefore the sequence 
$0\to \F \to \oplus^p{\mathcal O}_{\PP^n} \to \oplus^p{\mathcal O}_{\PP^n}/\F\to 0$
restricts to $H$ to give a short exact sequence
$0\to \F_H \to \oplus^p{\mathcal O}_H \to \oplus^p{\mathcal O}_H/\F_H \to 0$.  
This shows that $\F_H$ is isomorphic to a subsheaf of 
$\oplus ^p{\mathcal O}_{\PP^{n-1}_k}$ (under an identification of
$H$ with $\PP^{n-1}_k$),
which is a basic step needed for our inductive argument.

Note that $\F$ is torsion free if non-zero, and so 
we have a short exact sequence $0\to \F(-1) \to \F \to \F_H \to 0$. 
From the associated cohomology sequence we get 
$\chi(\F_H(r)) = \chi(\F(r)) - \chi(\F(r-1))
= \sum_{i=0}^n a_i {r\choose  i} - \sum_{i=0}^n a_i {r-1\choose i}
=  \sum_{i=0}^n a_i {r-1\choose i-1}= \sum_{j=0}^{n-1}b_j{r\choose j}$
where the coefficients $b_0,\ldots,b_{n-1}$ have  
expressions $b_j = g_j(a_0,\ldots,a_n)$ where the $g_j$ are polynomials
with integral coefficients independent of the field $k$ 
and the sheaf $\F$. (\textbf{Exercise}: Write down the $g_j$ explicitly.) 

By inductive hypothesis on $n-1$ 
there exists a polynomial $F_{p,n-1}(x_0,\ldots,x_{n-1})$
such that $\F_H$ is $m_0$-regular where
$m_0 = F_{p,n-1}(b_0,\ldots,b_{n-1})$.
Substituting $b_j = g_j(a_0,\ldots,a_n)$, we get 
$m_0 = G(a_0,\ldots,a_n)$, where $G$ is a polynomial with
integral coefficients independent of the field $k$ 
and the sheaf $\F$. 

For $m\ge m_0 -1$, we therefore get a long exact cohomology sequence
{\footnotesize
$$0\to H^0(\F(m-1))\to H^0(\F(m)) \stackrel{\nu_m}{\to} H^0(\F_H(m))
\to H^1(\F(m-1))\to H^1(\F(m)) \to 0 \to \ldots $$}
which for $i\ge 2$ gives isomorphisms
$H^i(\F(m-1))\stackrel{\sim}{\to} H^i(\F(m))$. As we have
$H^i(\F(m))=0$ for $m\gg 0$, these equalities show that 
$$H^i(\F(m)) = 0\mbox{ for all }i\ge 2\mbox{ and }m\ge m_0-2.$$
The surjections $H^1(\F(m-1))\to H^1(\F(m))$ show that
the function $h^1(\F(m))$
is a monotonically decreasing function of the variable $m$ for $m\ge m_0 - 2$. 
We will in fact show that for $m\ge m_0$, the function
$h^1(\F(m))$ is strictly decreasing  till its value reaches zero, 
which would imply that 
$$H^1(\F(m))=0 \mbox{ for }  m \ge m_0 + h^1(\F(m_0)).$$

Next we will put a suitable upper bound on $h^1(\F(m_0))$ to complete the
proof of the theorem. Note that 
$h^1(\F(m-1))\ge h^1(\F(m))$ for $m\ge m_0$, 
and moreover equality holds for some $m\ge m_0$ 
if and only if the restriction map 
$\nu_m : H^0(\F(m)) \to H^0(\F_H(m))$ is surjective.
As $\F_H$ is $m$-regular, it follows from
Remark \ref{surjectivity of restriction} that the restriction map 
$\nu_j : H^0(\F(j)) \to H^0(\F_H(j))$ is surjective for all
$j\ge m$, so $h^1(\F(j-1)) = h^1(\F(j))$ for all $j\ge m$.
As $h^1(\F(j))=0$ for $j\gg 0$, this establishes our claim that 
$h^1(\F(m))$ is strictly decreasing for
$m\ge m_0$ till its value reaches zero.

To put a bound on $h^1(\F(m_0))$, we use the fact that 
as $\F \subset \oplus^p{\mathcal O}_{\PP^n}$ we must have 
$h^0(\F(r)) \le p h^0({\mathcal O}_{\PP^n}(r)) = p{n + r\choose  n}$. 
From the already established fact that 
$h^i(\F(m)) = 0$ for all $i\ge 2$ and $m\ge m_0-2$, we now get
\begin{eqnarray*}
h^1(\F(m_0))&=& h^0(\F(m_0)) - \chi(\F(m_0)) \\
            &\le&p{n + m_0\choose  n}-\sum_{i=0}^n a_i {m_0\choose i}\\
            &=& P(a_0,\ldots a_n)
\end{eqnarray*}
where $P(a_0,\ldots, a_n)$ is a polynomial expression in 
$a_0,\ldots, a_n$, obtained by substituting 
$m_0 = G(a_0,\ldots, a_n)$ in the second line of the above (in)equalities.
Therefore, the coefficients of the corresponding polynomial
$P(x_0,\ldots, x_n)$ are again independent of the field $k$ 
and the sheaf $\F$. Note moreover that as
$h^1(\F(m_0))\ge 0$, we must have $P(a_0,\ldots, a_n) \ge 0$.

Substituting in an earlier expression, we get
$$H^1(\F(m))=0 \mbox{ for }  m \ge G(a_0,\ldots, a_n) + P(a_0,\ldots, a_n)$$
Taking $F_{p,n}(x_0,\ldots, x_n)$ to be $G(x_0,\ldots, x_n) + 
P(x_0,\ldots, x_n)$, and noting the fact that 
$P(a_0,\ldots, a_n) \ge 0$, we see that 
$\F$ is $F_{p,n}(a_0,\ldots, a_n)$-regular.
This completes the proof of the theorem. \qed

\medskip

\textbf{Exercise } Write down such polynomials $F_{p,n}$.

\section{Semi-Continuity and Base-Change}

\centerline{\large\bf Base-change without Flatness}

\medskip

The following lemma on base-change does not need any flatness hypothesis.
The price paid is that the integer $r_0$ may depend on $\phi$.

\begin{lemma}\label{base change without flatness}
Let $\phi : T\to S$ be a morphism of noetherian schemes, 
let $\F$ a coherent sheaf on $\PP^n_S$, and 
let $\F_T$ denote the pull-back of
$\F$ under the induced morphism $\PP^n_T\to \PP^n_S$. Let
$\pi_S : \PP^n_S \to S$ and $\pi_T : \PP^n_T \to T$
denote the projections. Then there exists an integer 
$r_0$ such that the base-change homomorphism 
$$\phi^* {\pi_S}_*\, \F(r) \to {\pi_T}_* \, \F_T(r)$$ 
is an isomorphism for all $r\ge r_0$.
\end{lemma}

\proof As base-change holds for open embeddings, using a finite affine open
cover $U_i$ of $S$ and a finite affine open cover $V_{i,j}$ of each
$\phi^{-1}(U_i)$ (which is possible by noetherian hypothesis), 
it is enough to consider the case where $S$ and
$T$ are affine.

Note that for all integers $i$, the base-change homomorphism  
$$\phi^*{\pi_S}_*{\mathcal O}_{\PP^n_S}(i) \to {\pi_T}_* {\mathcal O}_{\PP^n_T}(i)$$
is an isomorphism. Moreover, if $a$ and $b$ are any integers
and if $f : {\mathcal O}_{\PP^n_S}(a) \to {\mathcal O}_{\PP^n_S}(b)$ is any homomorphism
and $f_T : {\mathcal O}_{\PP^n_T}(a) \to {\mathcal O}_{\PP^n_T}(b)$ denotes its pull-back 
to $\PP^n_T$, then for all $i$ we have the following commutative diagram 
where the vertical maps are base-change isomorphisms.
$$\begin{array}{ccc}
\phi^*{\pi_S}_*{\mathcal O}_{\PP^n_S}(a+i) 
& \stackrel{\phi^*{\pi_S}_*f(i)}{\to} &
\phi^*{\pi_S}_*{\mathcal O}_{\PP^n_S}(b+i) \\
\dna & & \dna \\
{\pi_T}_* {\mathcal O}_{\PP^n_T}(a+i)
& \stackrel{{\pi_T}_*f_T(i)}{\to} &
{\pi_T}_* {\mathcal O}_{\PP^n_T}(b+i) 
\end{array}
$$

As $S$ is noetherian and affine, there exists an exact sequence
$$\oplus^p {\mathcal O}_{\PP^n_S}(a)\, \stackrel{u}{\to}  \, 
\oplus^q {\mathcal O}_{\PP^n_S}(b) \,\stackrel{v}{\to} \,\F \,\to \,0$$  
for some integers $a$, $b$, $p\ge 0$, $q\ge 0$. 
Its pull-back to $\PP^n_T$ is an exact sequence
$$\oplus^p {\mathcal O}_{\PP^n_T}(a)\, \stackrel{u_T}{\to}  \, 
\oplus^q {\mathcal O}_{\PP^n_S}(b) \,\stackrel{v_T}{\to} \,\F_T \,\to \,0$$ 
Let $\G = \ker(v)$ and let $\H = \ker(v_T)$. 
For any integer $r$, we get exact sequences
$${\pi_S}_* \oplus^p {\mathcal O}_{\PP^n_S}(a+r) \,\to \, 
  {\pi_S}_* \oplus^q {\mathcal O}_{\PP^n_S}(b+r) \,\to \,
  {\pi_S}_* \F(r)                       \,\to \,
  R^1{\pi_S}_* \G(r)$$ 
and 
$${\pi_T}_* \oplus^p {\mathcal O}_{\PP^n_T}(a+r)   \,\to \,
  {\pi_T}_* \oplus^q {\mathcal O}_{\PP^n_T}(b+r) \,\to \,
  {\pi_T}_* \F_T(r)                    \,\to \, 
  R^1{\pi_T}_* \H(r)$$ 
There exists an integer $r_0$ such that 
$R^1{\pi_S}_*\G(r) =0$ and 
$R^1{\pi_T}_*\H(r)=0$ for all $r\ge r_0$.
Hence for all $r\ge r_0$, we have exact sequences
$${\pi_S}_* \oplus^p {\mathcal O}_{\PP^n_S}(a+r) \,\stackrel{{\pi_S}_*u(r)}{\to}\, 
\oplus^q {\mathcal O}_{\PP^n_S}(b+r) \,\stackrel{{\pi_S}_*v(r)}{\to} \,
{\pi_S}_* \F(r) \,\to \, 0$$
and 
$${\pi_T}_* \oplus^p {\mathcal O}_{\PP^n_T}(a+r) \,\stackrel{{\pi_T}_*u_T(r)}{\to}\,
{\pi_T}_*   \oplus^q {\mathcal O}_{\PP^n_T}(b+r) \,\stackrel{{\pi_T}_*v_T(r)}{\to} \,
{\pi_T}_* \F_T(r) \,\to\, 0$$

Pulling back the second-last exact sequence under $\phi: T\to S$, 
we get the commutative diagram with exact rows
$$\begin{array}{ccccc}
\phi^*{\pi_S}_* \oplus^p {\mathcal O}_{\PP^n_S}(a+r) &
\stackrel{\phi^*{\pi_S}_*u}{\to} &
\phi^* \oplus^q {\mathcal O}_{\PP^n_S}(b+r) & 
\stackrel{\phi^*{\pi_S}_*v}{\to} &
\phi^*{\pi_S}_* \F(r) \to  0 \\
\dna &   & \dna & & \dna ~~~~~~~~ \\
{\pi_T}_* \oplus^p {\mathcal O}_{\PP^n_T}(a+r) & 
\stackrel{{\pi_T}_*u_T(r)}{\to} & 
{\pi_T}_*   \oplus^q {\mathcal O}_{\PP^n_T}(b+r) & 
\stackrel{{\pi_T}_*v_T(r)}{\to} & 
{\pi_T}_* \F_T(r) ~~ \to  0
\end{array}
$$
in which the first row is exact by the right-exactness of tensor product.
The vertical maps are base-change homomorphisms, the first two 
of which are isomorphisms for all $r$. 
Therefore by the five lemma, 
$\phi^*{\pi_S}_* \F(r) \to {\pi_T}_* \F_T(r)$ is an isomorphism
for all $r\ge r_0$.   
\qed

\medskip

{\footnotesize 
The following elementary proof of the above result is taken from Mumford [M]:
Let $M$ be the graded ${\mathcal O}_S$-module 
$\oplus_{m\in \Z}\, {\pi_S}_* \F(m)$, so that $\F = M^{\sim}$.
Let $\phi^*M$ be the graded ${\mathcal O}_T$-module which is the pull-back of
$M$. Then we have $\F_T = (\phi^*M)^{\sim}$.
On the other hand, let $N = \oplus_{m\in \Z}\, {\pi_T}_* \F_T(m)$,
so that we have $\F_T = N^{\sim}$. Therefore, in the category of
graded ${\mathcal O}_T[x_0,\ldots,x_n]$-modules, we get an induced equivalence
between $\phi^*M$ and $N$, which means the natural homomorphisms of
graded pieces $(\phi^*M)_m \to N_m$ are isomorphisms for all $m\gg 0$. 
$\square$
}

\bigskip

\centerline{\large\bf Flatness of $\F$ from Local Freeness of $\pi_*\F(r)$}

\medskip

\begin{lemma}\label{graded module is flat}
Let $S$ be a noetherian scheme and let $\F$ be a coherent sheaf on
$\PP^n_S$. Suppose that there exists some integer $N$ such that
for all $r \ge N$ the direct image $\pi_*\F(r)$ is locally free. 
Then $\F$ is flat over $S$.
\end{lemma}

\proof Consider the graded module $M = \oplus_{r\ge N} M_r$ over
${\mathcal O}_S$, where $M_r =  \pi_*\F(r)$. The sheaf $\F$ is isomorphic
to the sheaf $M^{\sim}$ on 
$\PP^n_S = {\bf Proj}_S\,{\mathcal O}_S[x_0,\ldots,x_n]$ made from the 
graded sheaf $M$ of ${\mathcal O}_S$-modules. 
As each $M_r$ is flat over ${\mathcal O}_S$, so is $M$. Therefore
for any $x_i$ the localisation $M_{x_i}$ is flat over ${\mathcal O}_S$. 
There is a grading on $M_{x_i}$, indexed by $\Z$, defined 
by putting $\deg(v_p/x_i^q) = p-q$ for $v_p \in M_p$ (this is well-defined).
Hence the component $(M_{x_i})_0$ of degree zero, 
being a direct summand of $M_{x_i}$, is again 
flat over ${\mathcal O}_S$. But by definition of $M^{\sim}$, this is just
$\Gamma(U_i,\F)$, where $U_i = {\bf Spec}_S\,{\mathcal O}_S[x_0/x_i,\ldots,x_n/x_i] 
\subset \PP^n_S$. As the $U_i$ form an open cover of $\PP^n_S$,
it follows that $\F$ is flat over ${\mathcal O}_S$. \qed

\medskip

\textbf{Exercise } Show that the converse of the above lemma holds: 
if $\F$ is flat over $S$ 
then  $\pi_*\F(r)$ is locally free for all sufficiently large 
$r$.

\bigskip

\centerline{\large\bf Grothendieck Complex for Semi-Continuity}

\medskip

The following is a very important basic result of Grothendieck,
and the complex $K^{\cdot}$ occurring in it is called the Grothendieck 
complex. 


\begin{theorem}\label{Grothendieck complex}
Let $\pi: X\to S$ be a proper morphism of noetherian
schemes where $S=\spec A$ is affine, 
and let $\F$ be a coherent ${\mathcal O}_X$-module
which is flat over ${\mathcal O}_S$. Then 
there exists a finite complex
$$0\to K^0 \to K^1 \to \ldots \to K^n\to 0$$
of finitely generated projective $A$-modules, together with
a functorial $A$-linear isomorphism 
$$H^p(X, \F \otimes_A M) 
\stackrel{\sim}{\to} H^p(K^{\cdot}\otimes_A M)$$
on the category of all $A$-modules $M$.
\end{theorem}

The above theorem is the foundation for all results about direct images and 
base-change for flat families of sheaves, such as Theorem
\ref{Base-change for flat sheaves}. 

As another consequence of the above theorem, we have the following.

\begin{theorem}\label{cokernel of dual} {\rm ([EGA] III 7.7.6) }
Let $S$ be a noetherian scheme and $\pi : X \to S$ a proper morphism.
Let $\F$ be a coherent sheaf on $X$ which is flat over $S$.
Then there exists a coherent sheaf $\Qq$ on $S$ together with
a functorial ${\mathcal O}_S$-linear isomorphism 
$$\theta_{\G} : \pi_*(\F \otimes_{{\mathcal O}_X} \pi^*\G) \to \hom_{{\mathcal O}_S}(\Qq,\G)$$ 
on the category of all quasi-coherent sheaves $\G$ on $S$.
By its universal property, 
the pair $(\Qq, \theta)$ is unique up to a unique isomorphism.
\end{theorem}

\proof If $S =\spec A$, then 
we can take $\Qq$ to be the coherent sheaf associated to the $A$-module
$Q$ which is the cokernel of the 
transpose $\pa^{\vee} : (K^1)^{\vee} \to (K^0)^{\vee}$ where
$\pa: K^0 \to K^1$ is the differential of any chosen Grothendieck
complex of $A$-modules 
$0\to K^0 \to K^1 \to \ldots \to K^n\to 0$ for the 
sheaf $\F$, whose existence is given by Theorem \ref{Grothendieck complex}. 
For any $A$-module $M$, 
the right-exact sequence 
$(K^1)^{\vee} \to (K^0)^{\vee} \to Q\to 0$ with $M$ 
gives on applying $Hom_A(-,M)$ a left-exact sequence 
$$0 \to Hom_A(Q, M) \to K^0\otimes_AM \to K^1\otimes_AM$$
Therefore by Theorem \ref{Grothendieck complex}, we have
an isomorphism 
$$\theta^A_M : H^0(X_A, \F_A \otimes_A M)  \to  Hom_A(Q, M)$$
Thus, the pair $(Q, \theta^A)$
satisfies the theorem when $S = \spec A$. More generally, we can cover
$S$ by affine open subschemes. Then on their overlaps,
the resulting pairs $(\Qq,\theta)$ glue together by their uniqueness.
\qed

A \textbf{linear scheme} $\VV \to S$ over a noetherian base scheme $S$ is
a scheme of the form $\Spec \sym_{{\mathcal O}_S} \Qq $ where $\Qq$ is a coherent
sheaf on $S$. This is naturally a group scheme. Linear schemes
generalise the notion of (geometric) vector bundles, which are the 
special case where $\Qq$ is locally free of constant rank. 

The \textbf{zero section} $\VV_0 \subset \VV$ of a linear scheme
$\VV = \Spec \sym_{{\mathcal O}_S}\Qq$ is the closed subscheme defined
by the ideal generated by $\Qq$. Note that the projection $\VV_0 \to S$
is an isomorphism, and $\VV_0$ is just 
the image of the zero section $0 : S \to \VV$ of the group-scheme.

\begin{theorem} {\rm ([EGA] III 7.7.8, 7.7.9) }
Let $S$ be a noetherian scheme and $\pi : X \to S$ a projective 
morphism. Let $\E$ and $\F$ be coherent sheaves on $X$.
Consider the set-valued contravariant functor $\Hom(\E,\F)$ on $S$-schemes,
which associates to any $T\to S$ the set of all ${\mathcal O}_{X_T}$-linear
homomorphisms $Hom_{X_T}(\E_T,\F_T)$ where $\E_T$ and $\F_T$ denote the 
pull-backs of $\E$ and $\F$ under the projection $X_T \to X$.
If $\F$ is flat over $S$, then the above functor is representable
by a linear scheme $\VV$ over $S$. 
\end{theorem}

\proof First note that if $\E$ is 
a locally free ${\mathcal O}_X$-module, then $\Hom(\E,\F)$ is
the functor $T\mapsto H^0(X_T, (\F\otimes_{{\mathcal O}_X}\E^{\vee})_T)$. 
The sheaf $\F\otimes_{{\mathcal O}_X}\E^{\vee}$ is again flat over $S$,
so we can apply Theorem \ref{cokernel of dual} to
get a coherent sheaf $\Qq$, such that we have
$\pi_*(\F \otimes_{{\mathcal O}_X}\E^{\vee}\otimes_{{\mathcal O}_X} \pi^*\G) = 
\hom_{{\mathcal O}_S}(\Qq,\G)$ for all quasi-coherent sheaves $\G$ on $S$.   
In particular, if $f: \spec R \to S$ is any morphism then
taking $\G = f_*{{\mathcal O}_R}$ we get
\begin{eqnarray*}
Mor_S(\spec R, \Spec \sym_{{\mathcal O}_S}\Qq ) & = & 
Hom_{{\mathcal O}_S\mod}(\Qq,  f_*{{\mathcal O}_R}) \\
& = & H^0(X, \F \otimes_{{\mathcal O}_X}\E^{\vee}\otimes_{{\mathcal O}_X} \pi^*  f_*{{\mathcal O}_R})\\
& = & H^0(X_R, (\F \otimes_{{\mathcal O}_X}\E^{\vee})_R)\\
&= &  Hom_{X_R}(\E_R, \F_R).
\end{eqnarray*}
This shows that $\VV = \Spec \sym_{{\mathcal O}_S} \Qq $
is the required linear scheme when $\E$ is locally free on $X$.
More generally for an arbitrary coherent $\E$, over
any affine open $U\subset S$ there exist vector bundles
$E_1$ and $E_0$ on $X_U$ and a right exact sequence
$E_1 \to E_0 \to \E \to 0$. (This is where we need projectivity
of $X\to S$. Instead, we could have assumed just properness together
with the condition that locally over $S$ we have such a resolution of
$\E$.) Then applying the above argument to the functors
$\Hom(E_1,\F)$ and $\Hom(E_0,\F)$, 
we get coherent sheaves $\Qq_1$ and $\Qq_0$ on
$U$, and from the natural transformation 
$\Hom(E_0,\F) \to \Hom(E_1,\F)$ induced by the homomorphism
$E_1\to E_0$, we get a homomorphism $\Qq_1 \to \Qq_0$. Let 
$\Qq_U$ be its cokernel, and put 
$\VV_U = \Spec \Sym_{{\mathcal O}_U} \Qq_U$. It follows 
from its definition (and the left exactness of Hom) 
that the scheme $\VV_U$ has
the desired universal property over $U$. Therefore all such
$\VV_U$, as $U$ varies over an affine open cover of $S$, patch together
to give the desired linear scheme $\VV$. 
(In sheaf terms, the sheaves $\Qq_U$ will patch together to give
a coherent sheaf $\Qq$ on $S$ with $\VV = \Spec \sym_{{\mathcal O}_S} \Qq $.) 
\hfill$\square$

\rem\label{vanishing scheme}
In particular, note that the zero section
$\VV_0 \subset \VV$ is where
the universal homomorphism vanishes. 
If $f \in Hom_{X_T}(\E_T,\F_T)$ defines a morphism 
$\varphi_f : T \to \VV$, then the inverse image $f^{-1}\VV_0$
is a closed subscheme $T'$ of $T$ with the universal property
that if $U\to T$ is any morphism of schemes such that
the pull-back of $f$ is zero, then $U\to T$ factors via $T'$.

\bigskip

\centerline{\large\bf Base-change for Flat Sheaves}

\medskip

The following is the main result of Grothendieck on base change for flat
families of sheaves, which is a consequence of 
Theorem \ref{Grothendieck complex}.

\begin{theorem}\label{Base-change for flat sheaves}
Let $\pi: X\to S$ be a proper morphism of noetherian
schemes, and let $\F$ be a coherent ${\mathcal O}_X$-module
which is flat over ${\mathcal O}_S$. 
Then the following statements hold:

(1) For any integer $i$ the function $s\mapsto \dim_{\kappa(s)}H^i(X_s,\F_s)$
is upper semi-continuous on $S$, 

(2) The function $s\mapsto \sum_i (-1)^i \dim_{\kappa(s)}H^i(X_s,\F_s)$
is locally constant on $S$.

(3) If for some integer $i$, there is some integer $d\ge 0$ such that
for all $s\in S$ we have $\dim_{\kappa(s)}H^i(X_s,\F_s) = d$, 
then $R^i\pi_*\F$ is locally free of rank $d$,
and $(R^{i-1}\pi_*\F)_s \to H^{i-1}(X_s,\F_s)$ is an isomorphism 
for all $s\in S$.

(4) If for some integer $i$ and point $s\in S$ the map
$(R^i\pi_*\F)_s \to H^i(X_s,\F_s)$ is surjective,
then there exists
an open subscheme $U\subset S$ containing $s$ such that
for any quasi-coherent ${\mathcal O}_U$-module $\G$ the natural homomorphism
$$(R^i {\pi_U}_*\F_{X_U})\otimes_{{\mathcal O}_U}\G\to 
R^i{\pi_U}_*(\F_{X_U}\otimes_{{\mathcal O}_{X_U}}{\pi_U}^*\G)$$
is an isomorphism, where $X_U = \pi^{-1}(U)$ and 
$\pi_U: X_U \to U$ is induced by $\pi$. In particular,
$(R^i\pi_*\F)_{s'} \to H^i(X_{s'},\F_{s'})$ is an isomorphism
for all $s'$ in $U$.

(5) If for some integer $i$ and point $s\in S$ the map
$(R^i\pi_*\F)_s \to H^i(X_s,\F_s)$ is surjective,
then the following conditions (a) and (b) are equivalent: 

~~~~(a) The map $(R^{i-1}\pi_*\F)_s \to H^{i-1}(X_s,\F_s)$ is surjective.

~~~~(b) The sheaf $R^i\pi_*\F$ is locally free in a 
        neighbourhood of $s$ in $S$.
\end{theorem}

See for example Hartshorne [H] Chapter III, Section 12 for a proof.
It is possible to replace the use of the formal function 
theorem in [H] (or the original argument in [EGA] based on 
completions) in proving the statement (4) above,
with an elementary argument based on applying
Nakayama lemma to the Grothendieck complex.


\section{Generic Flatness and Flattening Stratification}

\medskip

\centerline{\large\bf Lemma on Generic Flatness}

\medskip

\begin{lemma} \label{algebraic lemma for generic flatness}
Let $A$ be a noetherian domain, 
and $B$ a finite type $A$ algebra. Let $M$ be a finite $B$-module. 
Then there exists an $f\in A$, $f\ne 0$, such that 
the localisation $M_f$ is a free module over $A_f$.
\end{lemma}

\proof Over the quotient field $K$ of $A$, the $K$-algebra 
$B_K = K\otimes _AB$ is of finite type, and $M_K = K\otimes_AM$ is a 
finite module over $B_K$. Let $n$ be the dimension of the support
of $M_K$ over $\spec B_K$. We argue by induction on $n$, starting
with $n=-1$ which is the case when $M_K =0$. In this case, as
$K\otimes _AM = S^{-1}M$ where $S = A - \{ 0\}$, each $v\in M$
is annihilated by some non-zero element of $A$. Taking a finite
generating set, and a common multiple of corresponding annihilating 
elements, we see there exists an $f\ne 0$ in $A$ with $fM=0$.
Hence $M_f =0$, proving the lemma when $n=-1$.

Now let $n\ge 0$, and let the lemma be proved for smaller values.
As $B$ is noetherian and $M$ is assumed to be a finite $B$-module,
there exists a finite filtration
$$0=M_0 \subset \ldots  \subset M_r =M$$
where each $M_i$ is a $B$-submodule of $M$ such that for each $i\ge 1$
the quotient module $M_i/M_{i-1}$ is isomorphic to $B/\pp_i$ for 
some prime ideal $\pp_i$ in $B$. 

Note that if $0\to M'\to M\to M'' \to 0$ is a short exact
sequence of $B$-modules, and if $f'$ and $f''$ are non-zero elements of
$A$ such that $M'_{f'}$ and $M''_{f''}$ are free over
respectively $A_{f'}$ and $A_{f''}$, then $M_f$ is a free module over $A_f$
where $f=f'f''$. We will use this fact repeatedly.
Therefore it is enough to prove the result when 
$M$ is of the form $B/\pp$ for a prime ideal $\pp$ in $B$. 
This reduces us to the case where $B$ is a domain and $M=B$.

As by assumption $K\otimes_AB$ has dimension $n\ge 0$ (that is,
$K\otimes_AB$ is non-zero), the map $A\to B$ must be injective.
By Noether normalisation lemma, there exist elements 
$b_1,\ldots,b_n \in B$, 
such that $K\otimes_AB$ is finite over its subalgebra
$K[b_1,\ldots, b_n]$ and 
the elements $b_1,\ldots, b_n$ are algebraically independent
over $K$. (For simplicity of notation, we write $1\otimes b$ simply as $b$.)
If $g\ne 0$ in $A$ is chosen to be a `common denominator' for  
coefficients of equations of integral dependence 
satisfied by a finite set of
algebra generators for $K\otimes_AB$ over $K[b_1,\ldots, b_n]$, 
we see that $B_g$ is finite over $A_g[b_1,\ldots, b_n]$. 

Let $m$ be the generic rank of the finite module $B_g$ over the
domain $A_g[b_1,\ldots, b_n]$. Then we have a short exact
sequence of $A_g[b_1,\ldots, b_n]$-modules of the form
$$0\to A_g[b_1,\ldots, b_n]^{\oplus m} \to B_g \to T \to 0$$
where $T$ is a finite torsion module over $A_g[b_1,\ldots, b_n]$.
Therefore, the dimension of the support of $K\otimes_{A_g}T$
as a $K\otimes_{A_g}(B_g)$-module is strictly less than $n$. Hence by
induction on $n$ (applied to the data $A_g$, $B_g$, $T$), 
there exists some $h\ne 0$ in $A$ with $T_h$ free over $A_{gh}$.
Taking $f=gh$, the lemma follows from the 
above short exact sequence.
\hfill$\square$

The above theorem has the following consequence, which follows
by restricting attention to a non-empty affine open subscheme of $S$.

\begin{theorem}\label{generic flatness} 
Let $S$ be a noetherian and integral scheme.
Let $p: X\to S$ be a finite type morphism, and let $\F$ be a 
coherent sheaf of ${\mathcal O}_X$-modules. Then there exists a non-empty
open subscheme $U\subset S$ such that the restriction of $\F$ 
to $X_U = p^{-1}(U)$ is flat over ${\mathcal O}_U$.
\end{theorem}



\centerline{\large\bf Existence of Flattening Stratification}

\medskip

\begin{theorem}\label{flattening stratification} 
Let $S$ be a noetherian scheme, and let 
$\F$ be a coherent sheaf on the projective space $\PP^n_S$
over $S$. Then the set $I$ of Hilbert polynomials of
restrictions of $\F$ to fibers of $\PP^n_S\to S$
is a finite set. Moreover, for each $f\in I$ there exist 
a locally closed subscheme $S_f$ of $S$,  
such that the following conditions are satisfied.

\textbf{(i) Point-set: } The underlying set $|S_f|$ of $S_f$ consists of 
all points $s\in S$ where the Hilbert polynomial of the
restriction of $\F$ to $\PP^n_s$ is $f$. In particular, 
the subsets $|S_f|\subset |S|$ are disjoint, 
and their set-theoretic union is $|S|$.

\textbf{(ii) Universal property: } 
Let $S'=\coprod S_f$ be the coproduct of the $S_f$, and let
$i:S'\to S$ be the morphism induced by the inclusions
$S_f\hra S$. Then the sheaf $i^*(\F)$ on $\PP^n_{S'}$ 
is flat over $S'$. Moreover, $i:S'\to S$ has the universal property that 
for any morphism $\varphi :T\to S$ the pullback $\varphi^*(\F)$ on
$\PP^n_T$ is flat over $T$ if and only if $\varphi$ factors through
$i:S'\to S$.  The subscheme $S_f$ is 
uniquely determined by the polynomial $f$.

\textbf{(iii) Closure of strata: } 
Let the set $I$ of Hilbert polynomials be given a total 
ordering, defined by putting $f<g$ whenever $f(n) < g(n)$ for all
$n \gg  0$. Then the closure in $S$ of the subset $|S_f|$ 
is contained in the union of all $|S_g|$ where $f\le g$. 
\end{theorem}

\proof It is enough to prove the theorem for open subschemes of
$S$ which cover $S$, as
the resulting strata will then glue together by their universal property.

\textbf{Special case: } Let $n=0$, so that $\PP^n_S=S$. 
For any $s\in S$, the \textbf{fiber} $\F|_s$ of $\F$ over $s$
will mean the pull-back of $\F$ to the subscheme
$\spec \kappa(s)$, where $\kappa(s)$ is the residue field at $s$. 
(This is obtained by tensoring the stalk of $\F$ at $s$ with the
residue field at $s$, both regarded as ${\mathcal O}_{S,s}$-modules.)
The Hilbert polynomial of the 
restriction of $\F$ to the fiber over $s$ is the degree $0$ 
polynomial $e \in \Q[\lambda]$, where $e = \dim_{\kappa(s)}\F|_s$. 

By Nakayama lemma, any basis of $\F|_s$ prolongs to a neighbourhood $U$ of $s$ 
to give a set of generators for $F|_U$. Repeating this argument,
we see that there exists a smaller neighbourhood $V$ of $s$ in which
there is a right-exact sequence
$${\mathcal O}_V^{\oplus m} \stackrel{\psi}{\to} {\mathcal O}_V^{\oplus e} 
\stackrel{\phi}{\to} \F \to 0$$
Let $I_e\subset {\mathcal O}_V$ be the ideal sheaf
formed by the entries of the $e\times m$ matrix
$(\psi_{i,j})$ of the homomorphism 
${\mathcal O}_V^{\oplus m} \stackrel{\psi}{\to} {\mathcal O}_V^{\oplus e}$. 
Let $V_e$ be the closed subscheme of $V$ defined by $I_e$. 
For any morphism of schemes
$f: T\to V$, the pull-back sequence 
$${\mathcal O}_T^{\oplus m} \stackrel{f^*\psi}{\to} {\mathcal O}_T^{\oplus e} 
\stackrel{f^*\phi}{\to}f^* \F \to 0$$
is exact, by right-exactness of tensor products. 
Hence the pull-back $f^*\F$ is a locally free ${\mathcal O}_T$-module of
rank $e$ if and only if $f^*\psi =0$, that is, 
$f$ factors via the subscheme $V_e\hra V$ defined by the vanishing of all
entries $\psi_{i,j}$.
Thus we have proved assertions \textbf{(i)} and \textbf{(ii)} of the theorem.

As the rank of the matrix $(\psi_{i,j})$ is lower semi-continuous,
it follows that the function $e$ is upper semi-continuous, which
proves the assertion \textbf{(iii)} of the theorem, completing its proof
when $n=0$.

\textbf{General case: } We now allow the integer $n$ to be arbitrary.
The idea of the proof is as follows:  
We show the existence of a stratification
of $S$ which is a `g.c.d.' of the flattening
stratifications for direct images $\pi_*\F(i)$ for all $i \ge N$ for some
integer $N$ (where the flattening
stratifications for $\pi_*\F(i)$ exist by case $n=0$ which we have 
treated above). This is the desired flattening stratification of $\F$
over $S$, as follows from Lemma \ref{graded module is flat}.

As $S$ is noetherian, it is a finite union of irreducible components,
and these are closed in $S$.
Let $Y$ be an irreducible component of $S$, and let $U$ be the
non-empty open subset of $Y$ which consists of all points which do not 
lie on any other irreducible component of $S$. 
Let $U$ be given the reduced subscheme
structure. Note that this makes $U$ an integral scheme, which is a
locally closed subscheme of $S$. 
By Theorem \ref{generic flatness} on generic flatness, 
$U$ has a non-empty open subscheme $V$ such that
the restriction of $\F$ to $\PP^n_V$ is flat over ${\mathcal O}_V$.
Now repeating the argument with $S$ replaced by its reduced 
closed subscheme $S-V$,
it follows by noetherian induction on $S$ that there exist
finitely many reduced, locally closed, mutually disjoint subschemes
$V_i$ of $S$ such that set-theoretically 
$|S|$ is the union of the $|V_i|$ and 
the restriction of $\F$ to $\PP^n_{V_i}$ is flat over ${\mathcal O}_{V_i}$.
As each $V_i$ is a noetherian scheme, and as the Hilbert polynomials
are locally constant for a flat family of sheaves, it follows 
that only finitely many polynomials occur in $V_i$
in the family of Hilbert polynomials $P_s(m) = \chi(\PP^n_s,\F_s(m))$ 
as $s$ varies over points of $V_i$. This allows us to conclude 
the following:

\medskip

\textbf{(A) } Only finitely many distinct Hilbert polynomials 
$P_s(m) = \chi(\PP^n_s,\F_s(m))$ occur, as $s$ 
varies over all of $S$.

\medskip

By the semi-continuity theorem applied to the flat families 
$\F_{V_i} = \F|_{\PP^n_{V_i}}$ para\-met\-rised by the 
finitely many noetherian schemes 
$V_i$, we get the following:

\medskip
\textbf{(B) } There exists an integer $N_1$ such that 
$R^r\pi_*\F(m)=0$ for all $r\ge 1$ and $m\ge N_1$,
and moreover $H^r(\PP^n_s,\F_s(m)) = 0$
for all $s\in S$.

\medskip

For each $V_i$, by Lemma \ref{base change without flatness}
there exists an
integer $r_i \ge N_1$ with the property
that for any $m\ge r_i$ the base change homomorphism
$$(\pi_*\F(m))|_{V_i} \to {\pi_i}_*\F_{V_i}(m)$$ 
is an isomorphism, where $\F_{V_i}$ denotes the restriction of
$\F$ to $\PP^n_{V_i}$, and $\pi_i : \PP^n_{V_i}\to V_i$ the projection. 
As the higher cohomologies of all fibers (in particular, the first
cohomology) vanish by \textbf{(B)}, 
it follows by semi-continuity theory for the flat family $\F_{V_i}$
over $V_i$ that for any $s\in V_i$ 
the base change homomorphism
$$({\pi_i}_*\F_{V_i}(m))|_s \to  H^0(\PP^n_s,\F_s(m))$$
is an isomorphism for $m\ge r_i$. Taking $N$ to be the maximum of
all $r_i$ over the finitely many non-empty $V_i$, and composing the
above two base change isomorphisms, we get the following.

\medskip

\textbf{(C) } There exists an integer $N \ge N_1$ such that 
the base change homomorphism
$$(\pi_*\F(m))|_s \to  H^0(\PP^n_s,\F_s(m))$$
is an isomorphism for all $m\ge N$ and $s\in S$.

\medskip

\textbf{Note } We now forget the subschemes $V_i$ but retain the facts 
\textbf{(A)}, 
\textbf{(B)}, 
\textbf{(C)} 
which were proved using the $V_i$.

\medskip

Let $\pi : \PP^n_S \to S$ denote the projection. Consider the 
coherent sheaves $E_0, \ldots, E_n$ on $S$, defined by
$$E_i = \pi_*\F(N+i)\mbox{ for } i=0,\ldots,n.$$
By applying the special case of of the theorem (where the relative dimension
$n$ of $\PP^n_S$ is $0$) to the sheaf $E_0$ on $\PP^0_S=S$, we get 
a stratification $(W_{e_0})$ of $S$
indexed by integers $e_0$, such that for any morphism 
$f: T\to S$ the pull-back $f^*E_0$ is a locally free ${\mathcal O}_T$-module 
of rank $e_0$ if and only if $f$ factors via $W_{e_0} \hra S$. 
Next, for each stratum $W_{e_0}$, we take the 
flattening stratification $(W_{e_0,e_1})$ for $E_1|_{W_{e_0}}$, and so on.
Thus in $n+1$ steps, we obtain finitely many 
locally closed subschemes 
$$W_{e_0,\ldots,e_n} \subset S$$
such that for any morphism $f:T\to S$ the pull-back 
$f^*E_i$ for $i=0,\ldots,n$ is a locally free ${\mathcal O}_T$-modules of
of constant rank $e_i$ if and only if $f$ factors 
via $W_{e_0,\ldots,e_n} \hra S$.

For any integer $N$ and $n$ where $n\ge 0$, 
there is a bijection from the set of numerical polynomials
$f\in \Q[\lambda]$ of degree $\le n$ to the set $\Z^{n+1}$,
given by 
$$f \mapsto (e_0,\ldots,e_n)\mbox{ where }e_i = f(N+i).$$ 
Thus, each tuple $(e_0,\ldots,e_n) \in \Z^{n+1}$ can be uniquely 
replaced by a numerical polynomial
$f\in \Q[\lambda]$ of degree $\le n$, allowing us to re-designate 
$W_{e_0,\ldots,e_n} \subset S$ as $W_f \subset S$.

Note that at any point $s\in S$, 
by \textbf{(B)} we have
$H^r(\PP^n_s,\F_s(m)) = 0$ for all $r\ge 1$ and $m\ge N$. 
The polynomial $P_s(m) = \chi(\PP^n_s,\F_s(m))$ has degree $\le n$,
so it is determined by its $n+1$ values $P_s(N),\ldots, P_s(N+n)$.
This shows that at any point $s\in W_f$, the Hilbert polynomial
$P_s(m)$ equals $f$. The desired locally closed subscheme $S_f \subset S$, 
whose existence is asserted by the theorem, 
will turn out to be a certain closed subscheme $S_f \subset W_f$ whose 
underlying subset is all of $|W_f|$. The scheme structure of 
$S_f$ (which may in general differ from that of $W_f$) is defined as follows.  

For any $i\ge 0$ and $s\in S$, the base change homomorphism
$$(\pi_*\F(N+i))|_s \to  H^0(\PP^n_s,\F_s(N+i))$$
is an isomorphism by statement \textbf{(C)}.
Hence each $\pi_*\F(N+i)$ has fibers of constant rank $f(N+i)$
on the subscheme $W_f$. However, this does not mean $\pi_*\F(N+i)$ 
restricts to a locally constant sheaf of rank $f(N+i)$. But it means that 
$W_f$ has a closed subscheme $W_f^i$,
whose underlying set is all of $|W_f|$, such that
$\pi_*\F(N+i)$ is locally free of rank $f(N+i)$ when restricted to 
$W_f^{(i)}$, and moreover has the property that any base-change
$T\to S$ under which $\pi_*\F(N+i)$ pulls back to a 
locally free sheaf of rank $f(N+i)$ factors via $W_f^i$.
The scheme structure of $W_f^{(i)}$ is defined by a coherent ideal sheaf 
$I_i \subset {\mathcal O}_{W_f}$. 
Let $I\subset {\mathcal O}_{W_f}$ be the sum of the $I_i$ over $i\ge 0$. 
By noetherian condition, the increasing sequence 
$$I_0 \subset I_0 + I_1 \subset I_0 + I_1 + I_2 \subset \ldots $$
terminates in finitely many steps, showing $I$ is again a coherent
ideal sheaf. Let $S_f\subset W_f$ be the closed subscheme defined by the 
ideal sheaf $I$. Note therefore that $|S_f| = |W_f|$ 
and for all $i\ge 0$, the sheaf $\pi_*\F(N+i)$ is locally free of 
rank $f(N+i)$ when restricted to $S_f$.

It follows that from their definition that 
the $S_f$ satisfy property (i) of the theorem.

We now show that the morphism $\coprod_f S_f \,\to \,S$ 
indeed has the property (ii) of the theorem.
By Lemma \ref{base change without flatness}, 
there exists some $N' \ge N$ such that for all $i\ge N'$, 
the base-change
$(\pi_*\F(i))|_{S_f} \to (\pi_{S_f})_*\F_{S_f}(i)$ is an isomorphism for
each $S_f$. Therefore 
$\F_{S_f}$ is flat over $S_f$ by Lemma \ref{graded module is flat},
as the direct images $\pi_*\F(i)$ for all $i \ge N'$ are locally
free over $S_f$. Conversely, if $\phi : T\to S$ is a morphism such
that $\F_T$ is flat, then the Hilbert polynomial is locally 
constant over $T$. Let $T_f$ be the open and closed subscheme of $T$ 
where the Hilbert polynomial is $f$. Clearly, the set
map $|T_f| \to |S|$ factors via $|S_f|$. But as the direct images
${\pi_T}_*\F_T(i)$ are locally free of rank $f(i)$ on $T_f$, it
follows in fact that the schematic morphism $T_f \to S$ factors 
via $S_f$, proving the property (ii) of the theorem.

As by \textbf{(A)} only finitely many polynomials $f$ occur, there
exists some $p \ge N$ such that for any two polynomials $f$ and $g$ 
that occur, we have $f < g$ if and only if $f(p) < g(p)$. 
As $S_f$ is the flattening stratification for $\pi_*\F(p)$,
the property (iii) of the theorem follows from the corresponding property in
the case $n=0$, applied to the sheaf $\pi_*\F(p)$ on $S$.

This completes the proof of the theorem.  \hfill$\square$


\textbf{Exercise } What is the flattening stratification of $S$ for 
the coherent sheaf ${\mathcal O}_{S^{red}}$ on $S$, 
where $S^{red}$ is the underlying reduced scheme of $S$?

\section{Construction of Quot Schemes}

\medskip

\centerline{\large\bf Notions of Projectivity}

\medskip

Let $S$ be a noetherian scheme. Recall that as defined by Grothendieck, 
a morphism $X\to S$ is called a \textbf{projective morphism}
if there exists a coherent sheaf $E$ on $S$, together with a 
closed embedding of $X$ into $\PP(E) = \Proj \sym_{{\mathcal O}_S} E $ 
over $S$. Equivalently, $X\to S$ is projective when 
it is proper and there exists a relatively very ample
line bundle $L$ on $X$ over $S$. These conditions are related by
taking $L$ to be the restriction of ${\mathcal O}_{\PP(E)}(1)$ to $X$, or
in the reverse direction, taking $E$ to be the direct image of $L$ on $S$.
A morphism $X\to S$ is called \textbf{quasi-projective} if it factors as 
an open embedding $X\hra Y$ followed by a projective morphism $Y \to S$.

A stronger version of projectivity was introduced 
by Altman and Kleiman: a morphism
$X\to S$ of noetherian schemes is called \textbf{strongly projective} 
(respectively, \textbf{strongly quasi-projective}) if
there exists a vector bundle $E$ on $S$ together with a closed
embedding (respectively, a locally closed embedding)
$X\subset \PP(E)$ over $S$.

Finally, the strongest version of (quasi-)projectivity 
is as follows (used for example in
the textbook [H] by Hartshorne): a morphism
$X\to S$ of noetherian schemes is projective in the strongest sense 
if $X$ admits a (locally-)closed embedding into $\PP^n_S$ for some $n$.

Note that none of the three versions of projectivity
is local over the base $S$.

\medskip

\textbf{Exercises } (i) Gives examples to show that the above
three notions of projectivity are in general distinct. 

(ii) Show that if $X\to S$ is projective and flat, where $S$ is noetherian,
then $X\to S$ is strongly projective.

(iii) Note that if every coherent sheaf of ${\mathcal O}_S$-modules
is the quotient of a vector bundle, then projectivity over the base $S$ 
is equivalent to strong projectivity. If $S$ admits an ample line bundle 
(for example, if $S$ is quasi-projective over an affine base), 
then all three notions of projectivity over $S$ 
are equivalent to each other.

\bigskip

\newpage

\centerline{\large\bf Main Existence Theorems}

\medskip

Grothendieck's original theorem on Quot schemes, 
whose proof is outlined in [FGA] TDTE-IV, is the following.

\begin{theorem}\label{Grothendieck Quot} {\rm (Grothendieck)}  
Let $S$ be a noetherian scheme, $\pi: X\to S$ a projective morphism,
and $L$ a relatively very ample line bundle on $X$. Then for any 
coherent ${\mathcal O}_X$-module $E$ and any polynomial $\Phi \in \Q[\lambda]$,
the functor $\Quot_{E/X/S}^{\Phi,L}$ is representable
by a projective $S$-scheme $\quot_{E/X/S}^{\Phi,L}$.
\end{theorem}

Altman and Kleiman gave a complete and detailed proof of the existence of
Quot schemes in [A-K 2]. They could remove the noetherian hypothesis, 
by instead assuming strong (quasi-)projectivity of $X\to S$ 
together with an assumption about the nature of the coherent  
sheaf $E$, and deduce that the scheme $\quot_{E/X/S}^{\Phi,L}$ is then
strongly (quasi-)projective over $S$.

For simplicity, in these lecture notes we state and prove the result in [A-K 2]
in the noetherian context.

\begin{theorem}\label{Altman-Kleiman Quot} {\rm (Altman-Kleiman)}
Let $S$ be a noetherian scheme, $X$ a closed subscheme of 
$\PP(V)$ for some vector bundle $V$ on $S$, $L = {\mathcal O}_{\PP(V)}(1)|_X$,
$E$ a coherent quotient sheaf of $\pi^*(W)(\nu)$ where
$W$ is a vector bundle on $S$ and $\nu$ is an integer, and 
$\Phi \in \Q[\lambda]$. 
Then the functor $\Quot_{E/X/S}^{\Phi,L}$ is representable
by a scheme $\quot_{E/X/S}^{\Phi,L}$ which can be embedded over $S$ as a 
closed subscheme of $\PP(F)$ for some vector bundle $F$ on $S$.

The vector bundle $F$ can be chosen to be an exterior power of 
the tensor product of $W$ with a symmetric powers of $V$.
\end{theorem}

Taking both $V$ and $W$ to be trivial in the above, we get the following.

\begin{theorem}
If $S$ is a noetherian scheme, $X$ is a closed subscheme of 
$\PP^n_S$ for some $n\ge 0$, $L = {\mathcal O}_{\PP^n_S}(1)|_X$,
$E$ is a coherent quotient sheaf of $\oplus^p {\mathcal O}_X(\nu)$ 
for some integers $p\ge 0$ and $\nu$, and $\Phi \in \Q[\lambda]$,
then the the functor $\Quot_{E/X/S}^{\Phi,L}$ is representable
by a scheme $\quot_{E/X/S}^{\Phi,L}$ which can be embedded over $S$ as a 
closed subscheme of $\PP^r_S$ for some $r\ge 0$.
\end{theorem}

The rest of this section is devoted to proving 
Theorem \ref{Altman-Kleiman Quot}, with extra noetherian hypothesis.
At the end, we will remark on how the proof also
gives us the original version of Grothendieck.

\bigskip

\centerline{\large\bf Reduction to the case of 
$\Quot_{\pi^*W/\PP(V)/S}^{\Phi,L}$}

\medskip

It is enough to prove Theorem \ref{Altman-Kleiman Quot} in the special case
that $X = \PP(V)$ and
$E= \pi^*(W)$ where $V$ and $W$ are vector bundles on $S$, as
a consequence of the next lemma.

\begin{lemma}
(i) Let $\nu$ be any integer. 
Then tensoring by $L^{\nu}$ gives an isomorphism of
functors from $\Quot_{E/X/S}^{\Phi,L}$ to 
$\Quot_{E(\nu)/X/S}^{\Psi,L}$ where the polynomial 
$\Psi \in \Q[\lambda]$ is defined by $\Psi(\lambda) = \Phi(\lambda + \nu)$. 

(ii) Let $\phi : E \to G$ be a surjective homomorphism of 
coherent sheaves on $X$. Then the corresponding natural transformation
$\Quot_{G/X/S}^{\Phi,L} \to \Quot_{E/X/S}^{\Phi,L}$
is a closed embedding.
\end{lemma}

\proof The statement (i) is obvious. The statement (ii) just says that
given any locally noetherian scheme $T$ and a family
$\langle \F,q\rangle \in \Quot_{E/X/S}^{\Phi,L}(T)$, there exists
a closed subscheme $T' \subset T$ with the following universal property:
for any locally noetherian scheme $U$ and a morphism $f: U\to T$,
the pulled back homomorphism of ${\mathcal O}_{X_U}$-modules $q_U : E_U \to \F_U$
factors via the pulled back homomorphism $\phi_U : E_U \to G_U$
if and only if $U \to T$ factors via $T'\hra T$. 
This is satisfied by taking $T'$ to be the vanishing scheme
for the composite homomorphism $\ker(\phi) \hra  E \stackrel{q}{\to} \F$
of coherent sheaves on $X_T$ (see Remark \ref{vanishing scheme}), 
which makes sense here as both $\ker(\phi)$ and $\F$ are coherent 
on $X_T$ and $\F$ is flat over $T$. \qed

Therefore if $\Quot_{\pi^*W/\PP(V)/S}^{\Phi,L}$ is
representable, then for any coherent quotient
$E$ of $\pi^*W(\nu)|_X$, we can take $\quot_{E/X/S}^{\Phi,L}$ 
to be a closed subscheme of $\quot_{\pi^*W/\PP(V)/S}^{\Phi,L}$.

\bigskip

\centerline{\large\bf Use of $m$-Regularity}

\medskip

We consider the sheaf $E = \pi^*(W)$ on $X=\PP(V)$ where
$V$ is a vector bundle on $S$,
and take $L = {\mathcal O}_{\PP(V)}(1)$. 
For any field $k$ and a $k$-valued point 
$s$ of $S$, we have an isomorphism 
$\PP(V)_s \simeq \PP^n_k$ where $n=\rank(V)-1$,
and the restricted sheaf $E_s$ on $\PP(V)_s$ is isomorphic to 
$\oplus^p{\mathcal O}_{\PP(V)_s}$ where $p=\rank(W)$. 
It follows from Theorem \ref{Mumford}
that given any $\Phi \in \Q[\lambda]$, there exists an integer
$m$ which depends only on $\rank(V)$, $\rank(W)$ and $\Phi$,
such that for any field $k$ and a $k$-valued point $s$ of $S$,
the sheaf $E_s$ on $\PP(V)_s$ is $m$-regular, and for any  
coherent quotient $q: E_s \to \F$ on $\PP(V)_s$ with 
Hilbert polynomial $\Phi$, the sheaf $\F$ and the kernel sheaf
$\G\subset E_s$ of $q$ are both $m$-regular. 
In particular, it follows from the Castelnuovo Lemma \ref{Castelnuovo}
that for $r\ge m$, all cohomologies 
$H^i(X_s,E_s(r))$, $H^i(X_s,\F(r))$, and $H^i(X_s,\G(r))$ 
are zero for $i\ge 1$,
and $H^0(X_s,E_s(r))$, $H^0(X_s,\F(r))$, and $H^0(X_s,\G(r))$ 
are generated by their global sections.

From the above it follows by Theorem \ref{Base-change for flat sheaves}
that if $T$ is an $S$-scheme and
$q :E_T \to \F$ is a $T$-flat coherent quotient with Hilbert polynomial
$\Phi$, then we have the following, where $\G\subset E_T$ is
the kernel of $q$.

\textbf{(*) } The sheaves ${\pi_T}_*\G(r)$, ${\pi_T}_*E_T(r)$, 
${\pi_T}_*\F(r)$ are locally free of fixed ranks 
determined by the data $n$, $p$, $r$, and $\Phi$, 
the homomorphisms ${\pi_T}^*{\pi_T}_*(\G(r))\to \G(r)$, 
${\pi_T}^*{\pi_T}_*(E_T(r))\to E_T(r)$, 
${\pi_T}^*{\pi_T}_*(\F(r))\to \F(r)$ are 
surjective, and the higher direct images  
$R^i{\pi_T}_*\G(r)$, $R^i{\pi_T}_*E_T(r)$, 
$R^i{\pi_T}_*\F(r)$ are zero, for all $r\ge m$ and $i\ge 1$. 

\textbf{(**) }
In particular we have the following commutative diagram
of locally sheaves on $X_T$, in which both rows are exact, and all 
three vertical maps are surjective.

$$\begin{array}{ccccccc}
0  \to & {\pi_T}^*{\pi_T}_*(\G(r)) & \to 
& {\pi_T}^*{\pi_T}_* (E_T(r)) &
\to  &{\pi_T}^*{\pi_T}_* (\F(r))& \to   0\\ 
       & \dna                &     & \dna                &
     & \dna               &        \\
0  \to & \G(r)                  & \to & E(r) &
\to  & \F(r)                 & \to   0
\end{array}$$

\bigskip

\centerline{\large\bf Embedding Quot into Grassmannian}

\medskip

We now fix a positive integer $r$ such that $r\ge m$.
Note that the rank of ${\pi_T}_*\F(r)$ is $\Phi(r)$
and $\pi_*E(r) = W \otimes_{{\mathcal O}_S}\sym^r V$.
Therefore the surjective homomorphism
${\pi_T}_*E_T(r) \to {\pi_T}_*\F(r)$
defines an element of the set 
$\Grass( W \otimes_{{\mathcal O}_S}\sym^r V , \Phi(r))(T)$.
We thus get a morphism of functors 
$$\alpha : \Quot^{\Phi,L}_{E/X/S} \to 
\Grass( W \otimes_{{\mathcal O}_S}\sym^r V, \Phi(r))$$
It associates to $q: E_T \to \F$ the quotient
${\pi_T}_*(q(r)) : {\pi_T}_*E_T(r) \to {\pi_T}_*\F(r)$.
 
The above morphism $\alpha$ is injective because the quotient 
$q: E_T \to \F$ can be recovered from 
${\pi_T}_*(q(r)) : {\pi_T}_*E_T(r) \to {\pi_T}_*\F(r)$ 
as follows. 

If $G = \grass( W \otimes_{{\mathcal O}_S}\sym^r V, \Phi(r))$ with
projection $p_G : G \to S$, and
$u : {p_G}^*E \to \U$ denotes the universal quotient on $G$
with kernel $v : \K \to {p_G}^*E$,
then the homomorphism 
${\pi_T}^*{\pi_T}_*(\G(r))\to {\pi_T}^*{\pi_T}_*E_T(r)$ 
can be recovered 
from the morphism $T \to G$ as
the pull-back of $v : \K \to {p_G}^*E$.
Let $h$ be  the composite 
${\pi_T}^*{\pi_T}_*(\G(r))\to {\pi_T}^*{\pi_T}_*(E_T(r))\to E_T(r)$.
As a consequence of the properties of the diagram \textbf{(**)}, 
the following is a right exact sequence on $X_T$
$${\pi_T}^*{\pi_T}_*(\G(r))\stackrel{h}{\to} 
E_T(r) \stackrel{q(r)}{\to} \F \to 0$$
and so $q(r) : E_T(r) \to \F(r)$ can be recovered as the cokernel of $h$.
Finally, twisting by $-r$, we recover $q$, proving the desired 
injectivity of the morphism of functors 
$\alpha :\Quot^{\Phi,L}_{E/X/S} \to 
\Grass( W \otimes_{{\mathcal O}_S}\sym^r V, \Phi(r))$.

\bigskip

\centerline{\large\bf Use of Flattening Stratification}
\nopagebreak
\medskip
\nopagebreak
We will next prove that 
$\alpha :\Quot^{\Phi,L}_{E/X/S} \to 
\Grass( W \otimes_{{\mathcal O}_S}\sym^r V, \Phi(r))$
is relatively representable. In fact, we will show that
given any locally noetherian $S$-scheme $T$ and
a surjective homomorphism $f :  W_T \otimes_{{\mathcal O}_T}\sym^r V_T \to \J$
where $\J$ is a locally free ${\mathcal O}_T$-module of rank $\Phi(r)$,
there exists a locally closed subscheme $T'$ of 
$T$ with the following universal property \textbf{(F)} : 

\textbf{(F) } Given any 
locally noetherian $S$-scheme $Y$ and an $S$-morphism $\phi :Y\to T$,
let $f_Y$ be the pull-back of $f$, and  
let $\K_Y = \ker(f_Y) = \phi^* \ker(f)$. Let
$\pi_Y : X_Y \to Y$ be the projection, and let
$h : {\pi_Y}^* \K_Y \to E_Y$ be the composite map
$${\pi_Y}^* \K_Y \to  {\pi_Y}^* (W \otimes_{{\mathcal O}_S}\sym^r V)
= {\pi_Y}^*{\pi_Y}_*E_Y \to  E_Y$$ 
Let $q:   E_Y \to \F$ be the cokernel of $h$. 
Then $\F$ is flat over $Y$ with its Hilbert polynomial on all fibers
equal to $\Phi$ if and only if $\phi :Y\to T$
factors via $Y'\hra Y$.

The existence of such a locally closed subscheme $T'$ of 
$T$ is given by Theorem \ref{flattening stratification}, which 
shows that $T'$ is the stratum corresponding to 
Hilbert polynomial $\Phi$ for the 
flattening stratification over $T$ for the sheaf $\F$ on $X_T$.

When we take $T$ to be $\grass( W \otimes_{{\mathcal O}_S}\sym^r V, \Phi(r))$ 
with universal quotient $u : {p_G}^*E \to \U$, the corresponding
locally closed subscheme $T'$ represents
the functor $\Quot^{\Phi,L}_{E/X/S}$ by its construction.

Hence we have shown that $\Quot^{\Phi,L}_{E/X/S}$ is represented
by a locally closed subscheme of 
$\grass( W \otimes_{{\mathcal O}_S}\sym^r V, \Phi(r))$. As 
$\grass( W \otimes_{{\mathcal O}_S}\sym^r V, \Phi(r))$ embeds as a 
closed subscheme of $\PP(\bigwedge ^{\Phi(r)} W \otimes_{{\mathcal O}_S}\sym^r V)$,
we get a locally closed embedding of $S$-schemes 
$$\quot^{\Phi,L}_{E/X/S} \subset 
\PP(\bigwedge ^{\Phi(r)} (W \otimes_{{\mathcal O}_S}\sym^r V))$$
In particular, the morphism $\quot^{\Phi,L}_{E/X/S}\to S$ 
is separated and of finite type.

\bigskip

\centerline{\large\bf Valuative Criterion for Properness}

\medskip

The original reference for the following argument is EGA IV (2) 2.8.1. 

The functor $\Quot_{E/X/S}^{\Phi,L}$ satisfies the following valuative
criterion for properness over $S$: given any discrete valuation ring
$R$ over $S$ with quotient field $K$, the restriction map
$$\Quot_{E/X/S}^{\Phi,L}(\spec R) \to \Quot_{E/X/S}^{\Phi,L}(\spec K)$$
is bijective. This can be seen as follows.
Given any coherent quotient 
$q: E_K \to \F$ on $X_R$ which defines an element $\langle \F, q\rangle$
of $\Quot_{E/X/S}^{\Phi,L}(\spec K)$.
Let $\ov{\F}$ be the image of the composite homomorphism
$E_R \to j_*(E_K) \to j_*\F$ where $j : X_K \hra X_R$ is the open
inclusion. Let $\ov{q} : E_R\to \ov{\F}$ be the induced surjection.
Then we leave it to the reader to verify that 
$\langle \F, q\rangle$ is an element of $\Quot_{E/X/S}^{\Phi,L}(\spec R)$
which maps to $\langle \F, q\rangle$, and is the unique such element.
(Use the basic fact that being flat over a dvr 
is the same as being torsion-free.)

As $S$ is noetherian and as we have already shown that 
$\quot^{\Phi,L}_{E/X/S}\to S$ is of finite type,
it follows that $\quot^{\Phi,L}_{E/X/S}\to S$ is a proper morphism.
Therefore the embedding of
$\quot^{\Phi,L}_{E/X/S}$ into 
$\PP(\bigwedge ^{\Phi(r)} (W \otimes_{{\mathcal O}_S}\sym^r V))$
is a closed embedding.

This completes the proof of Theorem \ref{Altman-Kleiman Quot}.\hfill$\square$

\bigskip

\centerline{\large\bf The Version of Grothendieck}

\medskip

We now describe how to get Theorem \ref{Grothendieck Quot}
from the above proof. As $S$ is noetherian, we can find a common
$m$ such that given any field-valued point 
$s: \spec k\to S$ and a coherent quotient $q: E_s \to \F$ on $X_s$
with Hilbert polynomial $\Phi$,
the sheaves $E_s(r)$, $\F(r)$, $\G(r)$ (where $\G = \ker(q)$) are
generated by global sections and all their higher cohomologies vanish,
whenever $r\ge m$.
This follows from the theory of $m$-regularity, and semi-continuity. 

Because we have such a common $m$, we get as before an injective morphism
from the functor $\Quot_{E/X/S}^{\Phi,L}$ into the Grassmannian functor
$\Grass(\pi_*E(r), \Phi(r))$. The sheaf 
$\pi_*E(r)$ is coherent, but need not be the quotient of a 
vector bundle on $S$. Consequently, the scheme $\grass(\pi_*E(r), \Phi(r))$
is projective over the base, but not necessarily strongly projective.

Finally, the use of flattening stratification, which can be 
made over an affine open cover of $S$, gives 
a locally closed subscheme of $\grass(\pi_*E(r), \Phi(r))$ which represents
$\Quot_{E/X/S}^{\Phi,L}$, which is in fact a
closed subscheme by the valuative criterion. 
Thus, we get $\quot_{E/X/S}^{\Phi,L}$ as a projective scheme over $S$.

\section{Some Variants and Applications}

\medskip

\centerline{\large\bf Quot Scheme in Quasi-Projective case}

\medskip

\textbf{Exercise } 
Let $\pi : Z\to S$ be a proper morphism of noetherian schemes. 
Let $Y\subset Z$ be a closed subscheme, and let $\F$ be a coherent sheaf
on $Z$. Then there exists an open subscheme $S' \subset S$ 
with the universal property that a morphism $T\to S$ factors through
$S'$ if and only if the support of the pull-back
$\F_T$ on $Z_T = Z\times_ST$ is disjoint from
$Y_T = Y \times_ST$.

\medskip

\textbf{Exercise } As a consequence of the above, show the following:  
If $\pi : Z\to S$ is a proper morphism with $S$ noetherian, if
$X\subset Z$ is an open subscheme, and if
$E$ is a coherent sheaf on $Z$, then 
$\Quot_{E|_X/X/S}$ is an open subfunctor of $\Quot_{E/Z/S}$.

\medskip

With the above preparation, the construction of 
a quot scheme extends to the strongly quasi-projective case, to
give the following.

\begin{theorem}\label{Strongly quasi-projective quot}
{\rm (Altman and Kleiman)} 
Let $S$ be a noetherian scheme, $X$ a locally closed subscheme of 
$\PP(V)$ for some vector bundle $V$ on $S$, $L = {\mathcal O}_{\PP(V)}(1)|_X$,
$E$ a coherent quotient sheaf of $\pi^*(W)(\nu)|_X$ where
$W$ is a vector bundle on $S$ and $\nu$ is an integer, and 
$\Phi \in \Q[\lambda]$. 
Then the functor $\Quot_{E/X/S}^{\Phi,L}$ is representable
by a scheme $\quot_{E/X/S}^{\Phi,L}$ which can be embedded over $S$ as a 
locally closed subscheme of $\PP(F)$ for some vector bundle $F$ on $S$.
Moreover, the vector bundle $F$ can be chosen to be an exterior power of 
the tensor product of $W$ with a symmetric power of $V$.
\end{theorem}

\proof Let $\ov{X} \subset \PP(V)$ be the schematic closure of 
$X\subset \PP(V)$, and let $\ov{E}$ be the coherent sheaf on $\ov{X}$
defined as the image of the composite homomorphism
$$\pi^*(W)(\nu)|_{\ov{X}} \to j_*(\pi^*(W)(\nu)|_X) \to j_*E$$
Then we get a quotient $\pi^*(W)(\nu)|_{\ov{X}} \to \ov{E}$
which restricts on $X\subset \ov{X}$ to the given quotient
$\pi^*(W)(\nu)|_X\to E$. Therefore by the above exercise,
$\Quot_{E/X/S}$ is an open subfunctor of 
$\Quot_{\ov{E}/\ov{X}/S}$. Now the result follows from the 
Theorem \ref{Altman-Kleiman Quot}.\qed

\medskip

In order to extend Grothendieck's construction of a quot
scheme to the quasi-projective case, one first needs the following lemma
which is of independent interest.

\begin{lemma}\label{coherent prolongation}
Any coherent sheaf on an open subscheme of a noetherian scheme $S$ can be
prolonged to a coherent sheaf on all of $S$.
\end{lemma}

\proof First consider the case where $S = \spec A$ is affine, and
let $j: U\hra S$ denote the inclusion.
The quasi-coherent sheaf 
$j_*(\F)$ corresponds to the $A$-module $M = H^0(S,j_*(\F))$,
in the sense that $j_*(\F) = M^{\sim}$. 
Given any $u\in U$, there exist finitely many 
elements $e_1,\ldots,e_n \in M$
which generate the fiber $\F_u$ regarded as a vector space
over the residue field $\kappa(u)$. By Nakayama these elements 
will generate the stalks of $\F$ in an open neighbourhood of $u$ in $U$.
Therefore by the noetherian hypothesis, there exist finitely many
elements $e_1,\ldots,e_r \in M$ which generate the stalk of $\F$ 
at each point of $U$. If $N\subset M$ is the submodule generated by these 
elements, then $\G = N^{\sim}$ is a coherent prolongation of
$\F$ to $S = \spec A$, proving the result in the affine case.

In the general case, by the 
noetherian condition there exists a maximal coherent prolongation
$(U',\F')$ of $\F$. Then unless $U' = S$, we can obtain  
a further prolongation of $\F'$ by using the affine case. 
For, if $u\in S- U'$,
we can take an affine open subscheme $V$ containing $u$, and a 
coherent prolongation $\G'$ of $\F'|_{U'\bigcap V}$ to all of $V$,
and then glue together $\G'$ and $\F'$ along $U'\bigcap V$ to 
further prolong $\F'$ to $U'\bigcup V$, 
contradicting the maximality of $(U',\F')$. 
\qed

\begin{theorem}\label{Quasi-projective quot}{\rm (Grothendieck)}
Let $S$ be a noetherian scheme, $X$ a quasi-projective 
scheme over $S$, $L$ a line bundle on $X$ which is
relatively very ample over $S$, $E$ a quotient sheaf on $X$, 
and $\Phi \in \Q[\lambda]$. 
Then the functor $\Quot_{E/X/S}^{\Phi,L}$ is representable
by a scheme $\quot_{E/X/S}^{\Phi,L}$ which is
quasi-projective over $S$.
\end{theorem}

\proof By definition of quasi-projectivity of $X\to S$, note that
$X$ can be embedded over $S$ as a locally closed subscheme of 
$\PP(V)$ for some coherent sheaf $V$ on $S$, 
such that $L$ is isomorphic to ${\mathcal O}_{\PP(V)}(1)|_X$.
Let $\ov{X} \subset \PP(V)$ be the schematic closure of
$X$ in $\PP(V)$. This is a projective scheme over $S$, and 
$X$ is embedded as an open subscheme in it. 
By Lemma \ref{coherent prolongation} the
coherent sheaf $E$ has a coherent prolongation $\ov{E}$ to $\ov{X}$.
For any such prolongation $\ov{E}$, the functor
$\Quot_{E/X/S}$ is an open subfunctor of $\Quot_{\ov{E}/\ov{X}/S}$. 
Therefore the desired result now follows from Theorem \ref{Grothendieck Quot}.
\qed

\bigskip

\centerline{\large\bf Scheme of Morphisms}

\medskip

We recall the following basic facts about flatness.

\begin{lemma}\label{local criterion}
\textbf{(1) } Any finite-type flat morphism between noetherian schemes is open.

\textbf{(2) } Let $\pi : Y \to X$ be a finite-type morphism of  
noetherian schemes. Then all $y\in Y$ such that $\pi$ is
flat at $y$ (that is, ${\mathcal O}_{Y,y}$ is a flat ${\mathcal O}_{X,\pi(y)}$-module)
form an open subset of $Y$.

\textbf{(3) } Let $S$ be a noetherian scheme, and let $f:X \to S$ 
and $g: Y\to S$ be finite type flat morphisms.
Let $\pi : Y \to X$ be any morphism such that $g = f\circ \pi$. 
Let $y\in Y$, let $x=\pi(y)$, and let
$s = g(y) = f(x)$. If the restricted 
morphism $\pi_s : Y_s \to X_s$ between the 
fibers over $s$ is flat at $y\in Y_s$, then $\pi$ is flat at $y\in Y$. 
\end{lemma}

\proof See for example Altman and Kleiman [A-K 1] Chapter V. 
The statement \textbf{(3) } is a consequence of what is known as the 
\textbf{local criterion for flatness}. \qed

\begin{theorem}\label{isomorphism is open}
Let $S$ be a noetherian scheme, and let $f:X \to S$ 
and $g: Y\to S$ be proper flat morphisms.
Let $\pi : Y \to X$ be any projective morphism with
$g = f\circ \pi$. Then $S$ has open subschemes
$S_2\subset S_1 \subset S$ with the following universal properties:

\textbf{(a) } For any locally noetherian $S$-scheme $T$, the base change
$\pi_T : Y_T \to X_T$ is a flat morphism if and only if the 
structure morphism $T\to S$ factors via $S_1$. (This does not need
$\pi$ to be projective.)

\textbf{(b) } For any locally noetherian $S$-scheme $T$, the base change
$\pi_T : Y_T \to X_T$ is an isomorphism if and only if the 
structure morphism $T\to S$ factors via $S_2$.
\end{theorem}

\proof \textbf{(a) } By Lemma \ref{local criterion}.(2), all
$y\in Y$ such that $\pi$ is flat at $y$ form an
open subset $Y'\subset Y$. Then $S_1 = S - g(Y - Y')$ 
is an open subset of $S$ as $g$ is proper. 
We give $S_1$ the open subscheme structure induced from $S$.
It follows from the local criterion of flatness 
(Lemma  \ref{local criterion}.(3))
that $S_1$ exactly consists of all $s\in S$ such that 
the restricted morphism $\pi_s : Y_s \to X_s$ between the 
fibers over $s$ is flat. Therefore again by the 
local criterion of flatness, $S_1$ has the desired universal property.

\textbf{(b) } Let $\pi_1 : Y_1 \to X_1$ be the pull-back of
$\pi$ under the inclusion $S_1 \hra S$.
Let $L$ be a relatively very ample 
line bundle for the projective morphism $\pi_1 : Y_1\to X_1$. 
Then by noetherianness there exists an integer $m\ge 1$ such that
${\pi_1}_*L^m$ is generated by its global sections and $R^i{\pi_1}_*L^m =0$
for all $i \ge 1$. By flatness of $\pi_1$, it follows that 
${\pi_1}_*L^m$ is a locally free sheaf. Let $U\subset X_1$ be the
open subschemes such that 
${\pi_1}_*L^m$ is of rank $1$ on $U$. 
Finally, let $S_2 = S_1 - f(X_1 - U)$, which is open as
$f$ is proper. 
We give $S_2$ the induced open subscheme structure, 
and leave it to the reader to verify that 
it indeed has the required universal property \textbf{(b) }. \qed

\medskip

If $X$ and $Y$ are schemes over a base $S$, then for any
$S$-scheme $T$, 
an \textbf{$S$-morphism from $X$ to $Y$ parametrised by $T$} 
will mean a $T$-morphism from $X\times_ST$ to $Y\times_ST$.
The set of all such will be denoted by $\Mor_S(X,Y)(T)$.
The association $T\mapsto  \Mor_S(X,Y)(T)$
defines a contravariant functor $\Mor_S(X,Y)$ from $S$-schemes to $Sets$.

\medskip

\textbf{Exercise } Let $k$ be a field, let $S=\spec k[[t]]$, 
$X= \spec k = \spec (k[[t]]/(t))$, and let $Y = \PP^1_S$.
Is $\Mor_S(X,Y)$ representable?

\medskip

\begin{theorem}
Let $S$ be a noetherian scheme, let $X$ be a projective scheme over $S$,
and let $Y$ be quasi-projective scheme
over $S$. Assume moreover that $X$ is flat over $S$. Then the
functor $\Mor_S(X,Y)$ is representable by an open
subscheme $\mor_S(X,Y)$ of $\hilb_{X\times_SY/S}$.
\end{theorem}

\proof We can associate
to each morphism $f : X_T \to Y_T$ (where $T$ is a scheme over $S$)
its graph $\Gamma_T(f) \subset (X\times_SY)_T$, which is
closed in $(X\times_SY)_T$ by separatedness of $Y\to S$.
We regard $\Gamma_T(f)$ as a closed subscheme of $(X\times_SY)_T$ 
which is isomorphic to $X$ under the graph morphism
$(\id_X,f): X \to  \Gamma_T(f)$ and the projection $\Gamma_T(f) \to X$,
which are inverses to each other.
As $X$ is proper and flat over $S$, so is $\Gamma_T(f)$, therefore
this defines a set-map $\Gamma_T : \Mor_S(X,Y)(T) \to 
\Hilb_{X\times_SY/S}(T)$ which is functorial in $T$, so
we obtain a morphism of functors 
$$\Gamma : \Mor_S(X,Y) \to \Hilb_{X\times_SY/S}$$

Given any element of $\Hilb_{X\times_SY/S}(T)$, represented by a
family $Z \subset (X\times_SY)_T$, it follows by applying 
Theorem \ref{isomorphism is open}.(b) to
the projection $Z\to X$ that $T$ has an open subscheme
$T'$ with the following universal property: for any
base-change $U \to T$, the pull-back $Z_U\subset (X\times_SY)_U$
maps isomorphically on to $X_U$ under the projection
$p : (X\times_SY)_U\to X_U$ if and only if $U\to T$ factors via
$T'$. Note therefore that over $T'$, 
the scheme $Z_{T'}$ will be the graph of a uniquely determined 
morphism $X_{T'} \to Y_{T'}$. 

This shows that the morphism of functors 
$\Gamma : \Mor_S(X,Y) \to \Hilb_{X\times_SY/S}$
is a representable morphism which is an open embedding. 
Therefore a representing scheme
$\mor_S(X,Y)$ for $\Mor_S(X,Y)$ exists as an open
subscheme of $\hilb_{X\times_SY/S}$. \qed

\textbf{Exercise } Let $S$ be a noetherian scheme and
$X \to S$ a flat projective morphism. 
Consider the set-valued contravariant functor $\Aut_{X/S}$ 
on locally noetherian $S$-schemes, which associates
to any $T$ the set of all automorphisms of $X_T$ over $T$.
Show that this functor is representable by an open subscheme of
$\mor_S(X,X)$.

\medskip

\textbf{Exercise } Let $S$ be a noetherian scheme 
and $\pi : Z \to X$ a morphism of $S$-schemes,
where $X$ is proper over $S$ and $Z$ is quasi-projective over $S$.
Consider the set-valued contravariant functor 
$\Pi_{Z/X/S}$ on locally noetherian $S$-schemes, which associates
to any $T$ the set of all sections of $\pi_T : Z_T \to X_T$. 
Show that this functor is representable by an open subscheme of
$\hilb_{Z/S}$.

\medskip

\bigskip

\centerline{\large\bf Quotient by a Flat Projective Equivalence Relation}

\medskip

Let $X$ be a scheme over a base $S$.
A \textbf{schematic equivalence relation} on $X$ over $S$ will
mean an $S$-scheme $R$ together with a morphism
$f: R \to X\times_SX$ over $S$ such that for any $S$-scheme
$T$ the set map $f(T) : R(T) \to X(T)\times X(T)$
is injective and its image is the graph of an equivalence relation on 
$X(T)$. (Here, we denote by $Z(T)$ the set $Mor_S(T,Z)= h_Z(T)$ of
all $T$-valued points of $Z$, where $Z$ and $T$ are $S$-schemes.)

We will say that a morphism $q : X\to Q$ of $S$-schemes  
is a \textbf{quotient} for a 
schematic equivalence relation $f: R \to X\times_SX$ over $S$ if $q$
is a \textbf{co-equaliser} for the component morphisms
$f_1,f_2 : R  \stackrel{\to}{\scriptstyle\to} X$ of 
$f: R \to X\times_SX$. This means $q \circ f_1 = q\circ f_2$, 
and given any $S$-scheme $Z$ and an $S$-morphism $g : X \to Z$ such that
$g\circ f_1 = g\circ f_2$, there exists a unique $S$-morphism
$h : Q \to Z$ such that $g = h\circ q$. 
A schematic quotient $q: X\to Q$, when it exists, is
unique up to a unique isomorphism. 
\textbf{Exercise:} A schematic quotient, when it exists,
is necessarily an epimorphism in the category of $S$-schemes.

\medskip

\textbf{Caution } Even if $q: X\to Q$ is a schematic quotient for $R$, 
for a given $T$ the map $q(T) : X(T) \to Q(T)$ may not be a quotient
for $R(T)$ in the category of sets. The map $q(T)$ may fail to 
be surjective, and moreover it may identify two distinct equivalence classes.
\textbf{Exercise:} Give examples where such phenomena occur.

\medskip

We will say that the quotient $q : X\to Q$ is \textbf{effective}
if the induced morphism $(f_1,f_2) : R \to X\times_QX$ is an isomorphism
of $S$-schemes. 
In particular, it will ensure that distinct equivalence classes
do not get identified under $q(T) : X(T) \to Q(T)$. But 
$q(T)$ can still fail to be surjective, as in the
following example.

\medskip

\textbf{Exercise } Let 
$S = \spec \Z$, and let $X \subset \AA^n_{\Z}$ be 
the complement of the zero section of $\AA^n_{\Z}$. 
Note that for any ring $B$, an 
element of $X(\spec B)$ is a vector 
$u\in B^n$ such that at least one component of $u$ is invertible in
$B$. Show that $X\times_SX$ has a closed subscheme $R$
whose $B$-valued points for any ring $B$ are all pairs 
$(u,v)\in X(\spec B)\times X(\spec B)$ such that
there exists an invertible element $\lambda \in B^{\times}$
with $\lambda u = v$. Show that an effective quotient $q: X\to Q$ exists, 
where $Q = \PP^{n-1}_{\Z}$. However, show that 
$q$ does not admit a global section, and so 
$q(Q) : X(Q) \to Q(Q)$ is not surjective.

\medskip

The famous example by Hironaka 
(see Example 3.4.1 in Hartshorne [H] Appendix B)
of a non-projective smooth complete variety $X$ 
over $\C$ together with a schematic equivalence relation $R$ 
(for which the morphisms $f_i : R\to X$ are finite flat,
in fact, \'etale of degree $2$) shows that 
schematic quotients do not always exist. 
But under the powerful assumption of projectivity,
Grothendieck proved an existence result for quotients, to which we
devote the rest of this section.

We will need the following elementary lemma from Grothendieck's 
theory of faithfully flat descent (this is a special case of
[SGA 1] Expos\'e VIII Corollary 1.9). 
The reader can consult the lectures of Vistoli
[V] for an exposition of descent.

\begin{lemma}\label{descent of subschemes}
(1) Any faithfully flat quasi-compact morphism of schemes $f : X \to Y$ 
is an effective epimorphism, that is, $f$ is a co-equaliser for
the projections $p_1,p_2: X\times_YX   \stackrel{\to}{\scriptstyle\to} X$.

(2) Let $p : D\to H$ be a faithfully flat quasi-compact morphism.
Let $Z \subset D$ be a closed subscheme
such that 
$$p_1^{-1}Z = p_2^{-1}Z \subset D\times_HD$$ 
where $p_1,p_2 :  D\times_HD \stackrel{\to}{\scriptstyle\to} D$ 
are the projections, and $p_i^{-1}Z$ is the schematic 
inverse image of $Z$ under $p_i$. 
Then there exists a unique closed subscheme $Q$ of $H$
such that $Z = p^{-1}Q \subset D$. By base-change from $p : D\to H$, 
it follows that the induced morphism $p|_Z : Z\to Q$ 
is faithfully flat and quasi-compact.\hfill $\square$
\end{lemma}

The idea of using Hilbert schemes to make
quotients of flat projective equivalence relations is due to Grothendieck, 
who used it in his construction of a relative Picard scheme.
In set-theoretic terms, the idea is actually very simple:
Let $X$ be a set,
and $R\subset X\times X$ an equivalence relation on $X$. 
Let $H$ be the power set of $X$ (means the set of all subsets of $X$),
and let ${\varphi}: X \to H$ be the map which sends $x\in X$ to its
equivalence class $[x] \in H$. If $Q\subset H$ is the image of ${\varphi}$, 
then the induced map $q: X\to Q$ is the quotient of $X$ modulo $R$ in the 
category of sets. The scheme-theoretic analogue of the above is the
following theorem of Grothendieck, where the Hilbert scheme of $X$ 
plays the role of power set. The first detailed proof appeared in
Altman and Kleiman [A-K 2].

\begin{theorem}\label{Projective flat quotient}
Let $S$ be a noetherian scheme, and let
$X\to S$ be a quasi-projective morphism.
Let $f : R \to X\times_SX$ be a schematic 
equivalence relation on $X$ over $S$, such that 
the projections $f_1,f_2 : R \stackrel{\to}{\scriptstyle\to} X$  
are proper and flat. 
Then a schematic quotient $X\to Q$ exists over $S$. 
Moreover, $Q$ is quasi-projective over $S$,
the morphism $X \to Q$ is faithfully flat and projective, and 
the induced morphism $(f_1,f_2) : R\to X\times_QX$ is an isomorphism (the
quotient is effective). 
\end{theorem}

\proof (Following Altman and Kleiman [A-K 2]) 
The properness of $f_i$
together with separatedness of $X\to S$ implies properness of
$f : R \to X\times_SX$. 
Also, $f$ is functorially injective by definition
of a schematic equivalence relation. 
It follows that $f$ is a closed embedding, which
allows us to regard $R$ as a closed subscheme of $X\times_SX$ 
(\textbf{Exercise}: 
\textit{Any proper morphism of noetherian schemes, which is injective 
at the level of functor of points, is a closed embedding}). 
This defines an element $(R)$ of $\Hilb_{X/S}(X)$, as the projection
$p_2|_R = f_2$ is proper and flat.

By Theorem \ref{Quasi-projective quot}, 
there exists a scheme $\hilb_{X/S}$ which represents the
functor $\Hilb_{X/S}$. 
As the parameter scheme $X$ is noetherian and as the Hilbert polynomial
is locally constant, 
only finitely many polynomials $\Phi$ occur as Hilbert 
polynomials of fibers of $f_2: R \to X$ with respect to a
chosen relatively very ample line bundle $L$ on $X$ over $S$.
Let $H$ be the finite disjoint union of the corresponding 
open subschemes  
$\hilb_{X/S}^{\Phi,L}$ of $\hilb_{X/S}$.
Then $H$ is a quasi-projective scheme over $S$ as 
each $\hilb_{X/S}^{\Phi,L}$ is so by 
Theorem \ref{Quasi-projective quot}. 
The family $(R) \in  \Hilb_{X/S}(X)$
therefore defines a classifying morphism 
${\varphi}: X \to H$, with the property that $\varphi^*D = R$ 
where $D\subset X\times_SH$ denotes the restriction to $H$ 
of the universal family over $\hilb_{X/S}$,
and $\varphi^*D$ denotes $(\id_X\times\,{\varphi})^{-1} D$.
Also, note that the projection $p: D\to H$ is proper and flat.
If $X$ is non-empty then each fiber of $f_2 : R \to X$
is also non-empty as the diagonal $\Delta_X$ is contained in $R$,
and therefore the Hilbert polynomial of each fiber is non-zero.
Hence $p: D\to H$ is surjective, and so $p$ is faithfully flat.

For any $S$-scheme $T$, it follows from its definition that
a $T$-valued point of $D$ is a pair $(x,V)$ with
$x\in X(T)$ and $V \in H(T)$ such that
$$x\in V$$
where the notation ``$x\in V$'' more precisely means that
the graph morphism $(x, \id_T) : T \to X\times_ST$ 
factors via $V \subset X\times_ST$. With this notation,
we will establish the following crucial property:

\medskip

\textbf{(***)} \textit{For any $S$-scheme $T$, and 
$T$-valued points $x,y \in X(T)$,
the following equivalences hold:}
$(x,y) \in R(T) \Leftrightarrow x \in \varphi(y)
\Leftrightarrow \varphi(x) = \varphi(y) \in H(T)$.

\medskip

For this, note that for any $x,y\in X(T)$, the morphism 
$(x,y): T \to X\times_SX$  factors as the composite
$$T \stackrel{(x,\id_T)}{\to}X\times_S T 
\stackrel{\id_X \times\, y}{\to}  X\times_S X$$
As $\varphi^*D = R$ and $(\varphi\circ y)^*D = \varphi(y)$,
it follows that $y^*R =  \varphi(y)$, in other words,
the schematic inverse image of $R\subset X\times_S X$ under 
$\id_X\times\, y : X\times_S T \to X\times_S X$
is $\varphi(y)\subset X\times_S T$. Hence the above factorisation
of $(x,y) : T\to X\times_SX$ shows that 
\begin{eqnarray*}
&& (x,y) \in R(T) \\
& \Leftrightarrow 
& \mbox{~the~morphism~} (x,\id_T):T \to X\times_ST 
  \mbox{~factors~via~} \varphi(y)\subset X\times_ST \\
& \Leftrightarrow & x \in \varphi(y)
\end{eqnarray*}
Moreover, $x \in \varphi(x)$ as $\Delta_X \subset R$.
Therefore, if $\varphi(x) = \varphi(y)$ then
$x \in \varphi(y)$. 

It now only remains to prove that 
if $(x,y) \in R(T)$ then $\varphi(x) = \varphi(y)$,
that is, the subschemes 
$\varphi(x)$ and $\varphi(y)$ of $X\times_S T$ are identical.
Note that 
$\varphi(x) = (\varphi\circ x)^*D = x^*\varphi^*D = 
x^*R$ and similarly $\varphi(y) = y^*R$, therefore
we wish to show that $x^*R  = y^*R$. To show this
in terms of functor of points, for any $T$-scheme $u: U \to T$  
we just have to show that
$(x^*R)(U) = (y^*R)(U)$ as subsets of $(X\times_ST)(U)$.
As $x^*R$ is the inverse image of $R$ under 
$\id_X\times\, x : X\times_ST \to X\times_SX$, it follows that
a $U$-valued point of $x^*R$ is the same as
an element $z\in X(U)$
such that $(\id_X\times\, x)\circ (z,u) \in R(U)$.
But as $(\id_X\times\, x)\circ (z,u) = (z, x\circ u)$, it follows that
$$z \in (x^*R)(U) \Leftrightarrow (z, x\circ u) \in R(U)$$
As $R(U)$ is an equivalence relation on the set $X(U)$, and 
as by assumption $(x,y) \in R(T)$, we have 
$(x\circ u, y\circ u) \in R(U)$, and so by transitivity 
we have
$$z \in (x^*R)(U) \Leftrightarrow (z, x\circ u) \in R(U)
\Leftrightarrow (z, y\circ u) \in R(U) 
\Leftrightarrow z \in (y^*R)(U)$$
Hence the subschemes $x^*R$ and $y^*R$ of $X\times_ST$
have the same $U$-valued points
for any $T$-scheme $U$, and therefore $x^*R = y^*R$, as was to be shown. 
This completes the proof of the assertion \textbf{(***)}.

The graph morphism $(\id_X, {\varphi}) : X \to X\times_S H$
is a closed embedding as $H$ is separated over $S$.
As $\Delta_X \subset R$ and as 
${\varphi}^* D = R$, it follows that
$(\id_X, {\varphi})$ factors through $D\subset X\times_SH$.
Thus, we get a closed subscheme $\Gamma_{\varphi} \subset D$,
which is the isomorphic image of $X$ under 
$(\id_X, {\varphi})$. We wish to apply the 
Lemma \ref{descent of subschemes}.(2)
to the faithfully flat quasi-compact morphism 
$p: D \to H$ and the closed subscheme $Z = \Gamma_{\varphi} \subset D$.

Any $T$-valued point of 
$\Gamma_{\varphi}$ is a pair $(x, {\varphi}(x)) \in D(T)$ where  
$x\in X(T)$. Any $T$-valued point of $D\times_HD$ is a triple
$(x,y,V)$ where $x,y\in X(T)$ and $V \in H(T)$ 
such that $x, y \in V$.
Under the projections
$p_1,p_2 :  D\times_HD \stackrel{\to}{\scriptstyle\to} D$,
we have $p_1(x,y,V) = (x,V)$ and $p_2(x,y,V) = (y,V)$.
We now have 
\begin{eqnarray*}
&& p_1(x,y,V) \in \Gamma_{\varphi}(T)\\
& \Leftrightarrow & (x,V) \in \Gamma_{\varphi}(T) \mbox{ and }y\in V \\
& \Leftrightarrow & V = \varphi(x)\mbox{ and }y\in V \\
& \Leftrightarrow & y \in \varphi(x) = V \\
& \Leftrightarrow & \varphi(y) = \varphi(x) =V ~~(\mbox{by the property  
\textbf{(***)}).}
\end{eqnarray*}
Similarly, we have $p_2(x,y,V) \in \Gamma_{\varphi}(T)$ if and only if
$\varphi(y) = \varphi(x) =V$.
Therefore, $p_1(x,y,V) \in \Gamma_{\varphi}(T)$ if and only if 
$p_2(x,y,V) \in \Gamma_{\varphi}(T)$.
This holds for all $T$-valued points for all $S$-schemes $T$, and so 
$p_1^{-1}\Gamma_{\varphi} = p_2^{-1}\Gamma_{\varphi} \subset D\times_HD$.
Therefore by Lemma \ref{descent of subschemes}(2)  
there exists a unique closed subscheme $Q \subset H$ such that 
$\Gamma_{\varphi}$ is the pull-back of $Q$ under $D\to H$. 
Let $p : \Gamma_{\varphi} \to Q$ be the morphism induced by the restriction
to $\Gamma_{\varphi}$ of $p: D\to H$. 
Let $q: X\to Q$ be defined as the composite
$$X \stackrel{(\id_X,\varphi)}{\to} \Gamma_{\varphi} 
\stackrel{p}{\to} Q$$
Then note that the composite 
$X \stackrel{q}{\to} Q \hra H$ equals $\varphi$. 

We will now show that $q: X\to Q$ as defined above is the desired
quotient of $X$ by $R$, with the required properties.

\textbf{(i) Quasi-projectivity of $Q\to S$ }: This is satisfied as
$Q$ is closed in $H$ and $H$ is quasi-projective over $S$.

\textbf{(ii) Faithful flatness and projectivity of $q$ }: 
This follows 
by base change from the faithfully flat projective morphism $p: D\to H$,
as the following square is Cartesian. 
$$\begin{array}{ccc}
X & \stackrel{(\id_X,\varphi)}{\to} & D \\
{\scriptstyle q} \dna ~ & \square & ~ \dna {\scriptstyle p} \\
Q & \hra & H
\end{array}$$

\textbf{(iii) Exactness of $R \stackrel{\to}{\scriptstyle\to} X \to Q$
and the isomorphism $R \to X\times_QX$ }:
By \textbf{(***)}, for any $T$-valued points $x,y \in X(T)$, we have
$(x,y)\in R(T)$ if and only if $\varphi(x) = \varphi(y)$.
This shows that the composite $R \stackrel{f_1}{\to} X \stackrel{q}{\to} Q$
equals the composite $R \stackrel{f_2}{\to} X \stackrel{q}{\to} Q$, and
the induced morphism $(f_1,f_2) : R \to X\times_QX$ 
is an isomorphism, by showing these statements
hold at the level of functor of points. Under the isomorphism
$R \to X\times_QX$, the 
morphisms $f_1,f_2 : R  \stackrel{\to}{\scriptstyle\to} X$
become the projection morphisms $p_1,p_2 : X\times_QX
\stackrel{\to}{\scriptstyle\to} X$.
By Lemma \ref{descent of subschemes}, the
morphism $q: X \to Q$ is a co-equaliser for $p_1,p_2$, and so $q$
is a co-equaliser for $f_1,f_2$. 

This completes the proof of Theorem \ref{Projective flat quotient}.\qed

What Altman and Kleiman actually prove in [A-K 2] is 
a strongly projective form of the above theorem
(without a noetherian assumption), 
using the hypothesis of strong quasi-projectivity 
in the following places in the above proof: if $X\to S$ is
strongly quasi-projective then $H \to S$ will again be so
by Theorem \ref{Strongly quasi-projective quot}, and therefore
$Q$ will be strongly quasi-projective over $S$. Moreover, $D\to H$
will be strongly projective, and therefore by base-change
$X\to Q$ will be strongly projective. 
In the noetherian case, this gives us the following result.

\begin{theorem}
Let $S$ be a noetherian scheme, and let
$X\to S$ be a strongly quasi-projective morphism.
Let $f : R \to X\times_SX$ be a schematic 
equivalence relation on $X$ over $S$, such that 
the projections $f_1,f_2 : R \stackrel{\to}{\scriptstyle\to} X$  
are proper and flat. 
Then a schematic quotient $X\to Q$ exists over $S$. 
Moreover, the quotient is effective, 
the morphism $X \to Q$ is
faithfully flat and strongly projective, and $Q$ is strongly
quasi-projective over $S$. 
\end{theorem}

\bigskip
\bigskip

\bibliographystyle{amsalpha}

\end{document}